\documentclass[12pt,a4paper]{amsart}

\usepackage[T1]{fontenc}
\input xypic\usepackage{amsmath}
\usepackage{amssymb}

\usepackage{amssymb}
\usepackage{latexsym}
\usepackage{graphicx}
\usepackage{psfrag}

\newcommand{\newtextt}{}%{}

\def\char{\hbox{char}\,}

\newcommand{\R}{{\mathbb R}} 

\newcommand{\C}{{\mathbb C}} 
\newcommand{\N}{{\mathbb N}} 
\newcommand{\cal}{\mathcal } 
\renewcommand{\L}{{\mathcal L}} 
\newcommand{\proofbox}{\hfill{$\Box$}\medskip} 

\def\hat{\widehat}
\def\tilde{\widetilde}
\def \bfo {\begin {eqnarray*} }
\def \efo {\end {eqnarray*} }
\def \ba {\begin {eqnarray*} }
\def \ea {\end {eqnarray*} }
\def \beq {\begin {eqnarray}}
\def \eeq {\end {eqnarray}}

\def \dist {\hbox{dist}}
\def\diag{\hbox{diag }}
\def \det {\hbox{det}}

\def \e {\varepsilon}
\def \p {\partial}

\def\F{{\mathcal F}}

\def\Z{{\mathbb Z}}
\newtheorem{definition}{Definition}[section] 
\newtheorem{theorem}[definition]{Theorem} 
\newtheorem{lemma}[definition]{Lemma}

\begin{document}
\title[Rigidity of broken geodesic flow]{Rigidity of broken geodesic flow and inverse problems}
\date{}%Mar. 9, 2007}
%\title[Rigidity of broken geodesics]{Rigidity of broken geodesics and an 
%inverse problem for the radiative transfer equation. \\
%*DRAFT VERSION*}
\author{Yaroslav Kurylev}
\address{Yaroslav Kurylev, University of Loughborough
Department of Mathematical Sciences,  
Loughborough, Leicestershire, LE11 3TU, UK
}
\email{}
\author{Matti Lassas}
\address{Matti Lassas, Helsinki University of Technology,
Institute of Mathematics,
PO Box 1100, FIN--02015 TKK, Finland}
\email{}
\author{Gunther Uhlmann}
\address{
Gunther Uhlmann, University of Washington, 
Department of Mathematics,
Seattle, Washington 98195-4350, USA}
\email{}

\maketitle

{\bf Abstract.} Consider a broken geodesics $\alpha([0,l])$ on a compact Riemannian
manifold $(M,g)$ with boundary of dimension $n\geq 3$. 
The broken geodesics are unions of two geodesics with the
property that they have a common end point. 
Assume that for every broken geodesic $\alpha([0,l])$
starting at and ending to the boundary $\partial M$
we know the  starting point and direction $(\alpha(0),\alpha'(0))$,
the end point and direction $(\alpha(l),\alpha'(l))$, and the length $l$.
We show that this data determines uniquely, up to an isometry, the manifold $(M,g)$.
This result has applications in inverse problems on very heterogeneous
media for situations
where there are many scattering points in the medium, and arises
in several applications including geophysics and medical imaging. 
As
an example we consider the inverse problem for the
radiative transfer equation (or the linear transport equation)
 with a non-constant wave speed. Assuming
that the scattering kernel is everywhere positive, we
show that the boundary measurements determine the wave speed
inside the domain up to an isometry.
\smallskip

\noindent 
{\bf AMS classification:} 35J25, 58J45.

\smallskip
 
\noindent {\bf  Keywords:} Rigidity of Riemannian manifolds,
broken geodesics, inverse problems, radiative transfer.

\section{Introduction.}
\subsection{Main result}
Let us consider a compact Riemannian manifold $(M,g)$ with boundary 
of dimension $n\geq 3$. Let $SM$ denote its
unit tangent bundle. The classical boundary rigidity problem
is the following (see \cite{Cr2,C,G,Gu,LSU,Mu,O, PU1, SU1, SU2}): Assume
that we know the  
distances $d(x,y)$ of boundary points $x,y\in \p M$. Can we
determine the isometry type of the manifold $(M,g)$?  Michel 
\cite{M1,M2}
observed
that in the case of simple manifolds 
these distance functions also determine the values of the
bicharacteristic flow at boundary, the so-called scattering
relation or lens relation, that is,
\ba
& &{\cal L}=\{(x,\xi),(y,\zeta),t)\in SM\times SM\times \R:\ x,y\in \p M,\\ 
& &\quad \quad \quad \quad \quad (\gamma_{x,\xi}(t),\p_t\gamma_{x,\xi}(t))=(y,\zeta)\hbox{ for some }t\geq 0\}
\ea
where  $\gamma_{x,\xi}$ is the geodesic of $(M,g)$ 
that leaves from $x$ to
direction $\xi$ at $t=0$. In other words, ${\cal L}$ gives the
information when and where and in which direction a geodesic, sent from the boundary,
hits again the boundary. It was shown in \cite{Gu} under
some conditions (see
also \cite{Ale1,Ale2}) that the wave front set of
the scattering operator associated to the
wave equation for the Laplace-Beltrami operator of a smooth Riemannian
metric determines the scattering
relation. 
The natural conjecture is that for non-trapping manifolds the scattering
relation determines the isometry type of the manifold. For recent
progress on this problem see the survey papers \cite{PU2,SU3}.

In the case of a very heterogeneous media with many scattering points inside
the manifold one can obtain further information by looking at the propagation
of singularities of waves going through the manifold. This
is the broken scattering relation
or broken lens relation that we proceed to define. 

A broken geodesic (or, a once broken geodesic) 
is a path $\alpha=
\alpha_{x,\xi,z,\eta}(t)$, where $z=\gamma_{x,\xi}(s)\in M$ for some $s\geq 0$,
$\eta\in S_zM$, and
\ba%\label{eq: broken}
\alpha_{x,\xi,z,\eta}(t)=
\left\{ \begin{array}{ll}\gamma_{x,\xi}(t),&t<s,\\
\gamma_{z,\eta}(t-s),& t\geq s,
\end{array}\right.
\ea
(See Fig.\ 1.)
In Riemannian geometry
broken geodesics are considered  e.g.\
in the classical Ambrose theorem \cite{Ab},
which says that the parallel translations of the curvature 
tensor along broken geodesics determine uniquely a simply connected
Riemannian manifold.

{We denote by $\ell(\alpha_{x,\xi,z,\eta})\in
\R_+\cup\{\infty\}$ smallest $l>0$ such that
$\alpha_{x,\xi,z,\eta}(l)\in \p M$.}
Denote by $\nu$ the interior unit normal vector and by
\ba
& &\Omega_+=\{(x,\xi)\in SM:\ x\in \p M,\ (\xi,\nu)_g>0\},\\
& &\Omega_-=\{(x,\xi)\in SM:\ x\in \p M,\ (\xi,\nu)_g<0\}
\ea
the incoming and outgoing boundary directions respectively.

The boundary entering and exiting points of broken geodesics  define  the broken scattering relation,
\ba
R&=&\{(x,\xi),(y,\zeta),t)\in SM\times SM\times \R_+:(x,\xi)
\in \Omega_+,\ (y,\zeta)\in \Omega_-,
\\ 
& &\ \ \ \ \ t=\ell(\alpha_{x,\xi,z,\eta}), \ \ \hbox{and}
\\ 
& &\ \ \ \ \ 
(\alpha_{x,\xi,z,\eta}(t),\p_t\alpha_{x,\xi,z,\eta}(t))=(y,\zeta)\hbox{ for some }(z,\eta)\in SM \}.
\ea
Our main result is:

\begin{theorem}\label{main} Let $(M,g)$ be a compact Riemannian
manifold with a non-empty boundary of dimension $n\geq 3$.
Then {$\p M$ and the} broken scattering relation $R$ determines the isometry type 
of the manifold 
$(M,g)$ uniquely.
\end{theorem}

We remark that this result doesn't assume any a-priori condition on the
metric $g$
or the manifold $M$. The difficulty in 
proving the result lies in the possible complicated nature of the broken
geodesic
flow.
The proof of the theorem above and the other results stated in the
introduction are given in 
sections 2--3.

\begin{figure}[htbp]
\begin{center}
%\psfrag{1}{$(x_0,\xi_0)$}
%\psfrag{2}{$singsupp(u_0)$}
%\psfrag{4}{$singsupp(u_1)$}
%\psfrag{3}{$singsupp(u_2)$}
%\includegraphics[width=5cm]{kuva1.eps}% \label{pict2}%was 12 cm
\includegraphics[width=5cm]{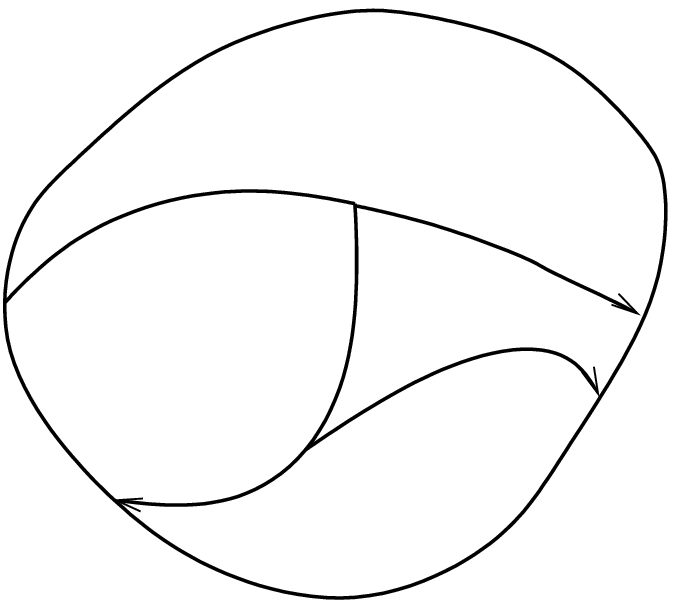}% \label{pict2}%was 12 cm
\psfrag{1}{$(x_0,\xi_0)$}
\psfrag{2}{$(x_1,\xi_1)$} 
\psfrag{3}{$s_1$} 
\psfrag{4}{$s_2$} 
\hspace{1cm}\includegraphics[width=5cm]{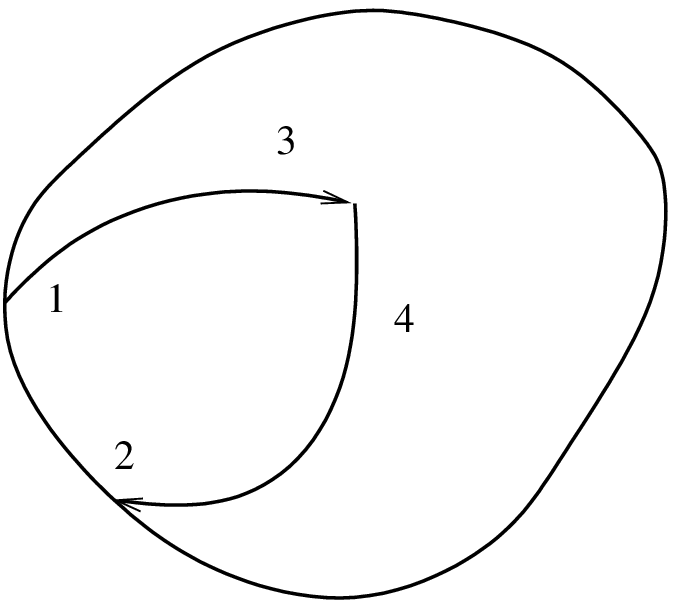}
\end{center}
\caption{Left: Propagation of singularities and multiple scattering for 
the radiative transfer equation.
Right: A broken geodesic corresponding the relation 
$((x_0,\xi_0),(x_1,\xi_1),t)\in R$ with $t=s_1+s_2$.}
\end{figure}

\subsection{Application: Radiative transfer equation}

As mentioned earlier the broken scattering relation can be determined by
probing with waves a very heterogeneous medium with many scattering points and observing at the boundary the effects.
The strongest singularities of the waves are the ones propagating through the medium
without any reflection and this determines the scattering
relation. The next stronger singularities correspond to the waves 
reflecting only once and this determines the broken scattering
relation at the boundary.  This type of situation
arises in geophysics due to the many discontinuities in the surface of the earth that act as reflectors
and in
optical tomography, a novel medical imaging technique that allows
one to reconstruct the spatial distribution of optical properties of tissues by probing them
by near-infra-red photons
\cite {Arr1,Arr2,Fe,He,HeS}.  This can be formulated as an inverse problem
for the radiative transfer equation and we consider this application
in more detail below. For previous mathematical analysis on the problem, see
e.g.\ \cite{Bal2,CS1,CS2,Hyv1,Hyv2, SU3}.

To avoid artificial difficulties on how to formulate the boundary 
value problem for the radiative transfer equation, we consider a non-compact
complete manifold  $(N,g)$ without boundary.
The inverse problem we study is to find the metric in a compact
subset $M$ with smooth boundary using  external measurements
made in the set  $U=N\setminus M$.

We say that the function $u(t,x,\xi)$ 
defined on
$(t,x,\xi)\in [0,\infty)\times SN$, is a solution
of the  radiative transfer equation on $N$
if 
\beq\label{a1}
& &(Hu)(t,x,\xi) +\sigma(x,\xi) u(t,x,\xi)-(Su)(t,x,\xi)=0,\\
& &u(t,x,\xi)|_{t=0}=w(x,\xi). \nonumber
\eeq
Here $H$ is the bicharacteristic flow on the tangent bundle $TN$,
\ba
Hu(t,x,\xi)=\frac {\p u}{\p t}-\xi^i\frac {\p u}{\p x^i}-\xi^i\xi^j 
\Gamma_{ij}^k(x)\frac {\p u}{\p \xi^k},
\ea
where $(x^1,\dots,x^n,\xi_1,\dots,\xi_n)$ denotes local coordinates
on the tangent bundle $TN$ corresponding
to local coordinates $(x^1,\dots,x^n)$ of $M$ and $\xi^j=g^{jk}\xi_k$.
The operator $S$, called the scattering
operator, is 
\ba%\label{scat. op}
Su (t,x,\xi)=c_n^{-1}\int_{S_xN}K(x,\xi,\xi')u(t,x,\xi')\,dV_g(\xi').
\ea
Here  $K\in C^\infty(SN\otimes SN)$ is called the scattering kernel and
$c_n=\hbox{vol}(S^{n-1})$. Finally, the function
$\sigma\in C^\infty(SN)$ is called the attenuation function.
We denote the solution of (\ref{a1}) with the initial
value $w\in C^\infty(SN)$ by $u(t,x,\xi)=
 u^w(t,x,\xi)$.

For the results concerning the radiative transfer equation 
we need a few more definitions.
We say that the complete manifold $N$ is simple if 
for any $x,y\in N$ there is only one geodesic
connecting these points.
{We say that $M\subset N$ is strictly convex if all points
in $M$ can be connected with a geodesic segment lying in $M$
and the second fundamental form of $\p M$ is positive.}  
%(We note that the analysis below 
%could be generalized to the case where only
% $M$ is assumed to be simple).

We say that scattering kernel $K$ is positive in $M^{int}$ if
\ba%\label{positive}
K(x,\xi,\xi')>0,\quad \hbox{for all $x\in M^{int}$ and $\xi,\xi'\in S_xN$.}
\ea

Next we define the external measurements.
We assume that for any $w\in C^\infty_0(SN)$,
such that $w(x,\xi)=0$ for $x\in M$ we know 
solution $u^w(x,\xi,t)$ for $x\in U$. In other words,
we assume that we are given
the measurement map $A:C^\infty_0(SU)\to C^\infty(\R_+\times SU),$ 
\ba
Aw=u^w|_{\R_+\times SU}.
\ea
Note that the map $A$ gives us the geodesic flow in $U$
and thus it determines the metric $g_{ij}(x)$ for $x\in U$. Also, it
can be used to determine the absorption $\sigma|_U$.

\begin{theorem}\label{transport-main}
Let $N$ be a complete simple manifold,
$M\subset N$ a  compact {and strictly convex} set with smooth boundary.
Assume that
$K(x,\theta,\theta')$ vanish for $x\not \in M$,  that is,
$K\in C^\infty_0(SM\otimes SM)$ and that $K$ is positive in $M^{\rm int}$.

Moreover, assume that we are given the set $U=N\setminus M$
and the measurement map $A$. These data 
determine uniquely the broken scattering 
relation of the manifold
$(M,g)$.
\end{theorem}

\section{Proof of Theorem \ref{main}}

\subsection{Auxiliary Lemmata}

Let $(M,g)$ be a compact manifold with boundary, $\p M$.
In the following,
we use an auxiliary smooth closed compact $n$--manifold
$(\tilde M,\tilde g)$ that contains $(M,g)$.
 We continue to use notation  $\gamma_{x,\xi}(t), \,
(x, \xi) \in S\tilde M,$ for  the geodesics on $\tilde M$
with $\gamma_{x,\xi}(0)=x$ and $\gamma'_{x,\xi}(t)=\xi$.
All geodesics are parameterized by  the arclength. 
We denote by
$\dist_{\tilde M}(x,y)$ and  $ \dist(x,y)$ 
the distance functions on $\tilde{M}$ and $M$, respectively.
To simplify notations, we denote
\ba
(x_0,\xi_0)R_t (x_1,\xi_1)\quad
\hbox{if and only if}\quad \bigg((x_0,\xi_0),(x_1,-\xi_1),t\bigg)\in R.
\ea

On $\tilde M$ and $M$, we will
use various critical distances along geodesics.
We  start with critical distances associated with the  Riemann exponential map,
$\exp_x,$
\ba
%\label{28.4.0}
\exp_x: T_x\tilde{M} \equiv S_x\tilde{M}
\times {\mathbb R}_+ \longrightarrow \tilde{M}, \quad
\exp_x(s\xi) = \gamma_{x,\xi}(s), 
\ea
$\xi \in S_x\tilde{M},\, s \in \R_+.$
The {\it cut locus distance} along $\gamma_{x,\xi}$, denoted
by $\tau_R(x,\xi)$,  is defined by 
\beq
\label{28.4.2}
\tau_R(x,\xi)=\max \{s>0:\, {\dist}_{\tilde M}(x,\gamma_{x,\xi}(s))=s\}.
\eeq
The cut locus distance 
$\tau_R(x,\xi), \,
(x,\xi) \in S\tilde{M}$ determines the 
injectivity radius $\hbox{inj}\,(M)$ of $\tilde{M}$,
\ba
%\label{24.5.2}
\hbox{inj}\,(M)= \min_{(x,\xi) \in S\tilde{M}} \tau_R (x,\xi).
\ea
We say that the set
\bfo
\omega_x = \{y \in \tilde M:\, y =\gamma_{x,\xi}(\tau_R(x,\xi)), \, \xi \in S_x\tilde M\},
\efo
is the {\it  cut locus} with respect to $x$. 
The cut locus
$\omega_x$ consists of two types of points.
We say that a point $y \in \omega_x$ is an 
{\it ordinary cut locus point} if there are $\xi, \eta \in S_x\tilde{M}$, $\eta\not=\xi$ 
with 
\ba
%\label{28.4.4}
\tau_R(x,\xi)= \tau_R(x,\eta), 
\quad \gamma_{x,\xi}(\tau_R(x,\xi)) =\gamma_{x,\eta}( \tau_R(x,\eta)) =y.
\ea
Consider now the differential of $\exp_x$ at ${s \xi }$
that is denoted by $d\exp_x|_{s \xi }$.
We say that a point 
$y=\gamma_{x,\xi}(s)$
 is a {\it conjugate point} along $\gamma_{x,\xi}$, if 
the differential 
 $d\exp_x|_{s \xi}:T_x\tilde M\to T_y\tilde M$ is degenerate.
This is equivalent to the existence
of   
a non-trivial Jacobi field $Y(t)$ along $\gamma=\gamma_{x,\xi}([0,s])$
with the Dirichlet boundary conditions $Y(0) = 0$ and $Y(s)=0$.
For $(x,\xi)\in S\tilde M$
we define the
 {\it conjugate distance} $\tau_c(x,\xi)\in \R_+\cup\{\infty\}$ to be
\ba
%\label{28.4.1}
\tau_c(x,\xi)=\inf \{s>0:\,d\exp_x|_{s \xi } \,\, \hbox{is not one-to-one}\}.
\ea
Each point $y\in \omega_x$ is an ordinary cut locus point, a first conjugate
point, or both.

\medskip

Next we discuss  critical distances associated with the boundary 
exponential map, $\exp_{\p M}$,
\ba
%\label{28.4.5}\quad\quad
\exp_{\p M}:\p M\times \R \longrightarrow \tilde M, \quad
\exp_{\p M}(z,s)=\gamma_{z,\nu}(s), \quad z \in \p M,
\ea
where $\nu = \nu (z)$ is the unit interior normal vector 
to $\p M$ at $z$. The pair $(z,s)$ defines the {\it boundary normal
coordinates} in $\tilde M$ near $\p M$.

The {\it boundary cut locus distance}, $\tau_b(z)$, $z\in \p M$ is given by  
\beq
\label{28.4.8}
\tau_b(z)=\max\{s>0:\ \dist(\gamma_{z,\nu}(s),\p M)=s\}.
\eeq
The set of the corresponding points
$
y= \gamma_{z,\nu}(\tau_b(z))$ is called the
{\it boundary cut locus}, 
\bfo
\omega_{\p M} = \{y \in  M:\, y =\gamma_{z,\nu}(\tau_b(z)), \, z \in \p M\}.
\efo
The boundary cut locus
consists of two types of points.
We say that a point $y \in \omega_{\p M}$ is an {\it ordinary boundary cut locus point}
 if there are 
$z, w \in \p M$, $z\not=w$ 
with 
\ba
%\label{28.4.10}
\tau_b(z) =\tau_b(w), 
\quad \gamma_{z,\nu(z)}(\tau_b(z)) =\gamma_{w,\nu(w)}( \tau_b(w)) =y.
\ea
Also, we say that a point  $y=\gamma_{z,\nu(z)}(\tau_b(z))\in \omega_x$ 
is a {\it 
focal point} if the differential,
$d \exp_{\p M}|_{(z,\tau_b(z))}:T_z\p M\times \R\to T_y\tilde M$ is 
degenerate. 
Equivalently, $t$ is a focal point if 
there is a non-trivial Jacobi field $Y(t)$ 
along $\gamma_{z, \nu}([0,s])$ with
 $Y(s)=0$ and $Y'(0) = WY(0)$,
where $W$ is the Weingarten map of $\p M$ at $z$. For $z\in \p M$,
we define the {\it focal distance},
$\tau_f(z)$ to be 
\ba
\label{28.4.6}
\tau_f(z)=\inf \{s>0:\, d \exp_{\p M}|_{(z,s)}  \,\, \hbox{is not one-to-one}\}.
\ea

Note that  $y\in \omega_{\p M}$ is an ordinary boundary cut locus 
point,
a first focal point, or both. Also, the functions
$\tau_R$, $\tau_c$, $\tau_b$, and $\tau_f$ are continuous,
e.g. \cite{Klingenberg}.

Comparing Jacobi fields $Y(s)$ along the geodesic
$\gamma_{z, \nu}([0,s])$ 
with the Dirichlet condition $Y(0)=0$ and the 
Robin condition  $Y'(0) = WY(0)$,
we see that $\tau_f(z)<\tau_c(z,\nu)$. Due
to the compactness of
 $\p M$ there is $c_0>0$ such that
 \ba
%\label{c_0-result}
\tau_c(z,\nu)\geq \tau_f(z)+c_0,\quad z\in \p M.
\ea

In a similar manner, we can show that 
$
\tau_R(z,\nu) > \tau_b(z),$ $z \in \p M.$ 
{Indeed, 
assume the opposite, i.e., $t=\tau_R(z, \nu) \leq \tau_b(z)$ 
for some $z \in \p M$. Denote 
$(y, \eta)=(\gamma_{z,\nu}(t),\,-\gamma'_{z,\nu}(t)).$
By duality, $\tau_R(y, \eta)=\tau_R(z, \nu)=t$. 
Let $\e >0$ and 
$x_{\e} =\gamma_{z,\nu}(-\e)=\gamma_{y, \eta}(t+\e).$
Then 
\ba
\dist_{\tilde M}\,(x_{\e},y) < t+\e \leq \tau_b(z)+\e
\ea
and there is $\eta_\e \in S_{x_\e}{\tilde M}$ with
$y= \gamma_{x_\e, \eta_\e}(\dist_{\tilde M}(x_{\e},y))$.
Denote by $t_\e>0$ the last time when $ \gamma_{x_\e, \eta_\e}(s)$
hits $\p M$.  If $\e$ is sufficiently small, we see by the
 short-cut arguments that
$\dist(y,\p M) < \tau_b(z)$. This contradicts the 
 definition  of $\tau_b$ in (\ref{28.4.8}).
} 

Due to the compactness of $\p M$, by making  $c_0>0$ smaller if
necessary,  
\beq
\label{24.5.7}
\tau_R(z,\nu) \geq \tau_b(z) +c_0, \quad z \in \p M.
\eeq

\begin{figure}[htbp]\label{same pointpicture3}
\begin{center}
\psfrag{2}{$\p M$}
\psfrag{1}{$z$} 
\psfrag{3}{$A$}
\psfrag{4}{$B$}
\psfrag{8}{$z_0$}
\psfrag{9}{$z$}
\includegraphics[height=5cm]{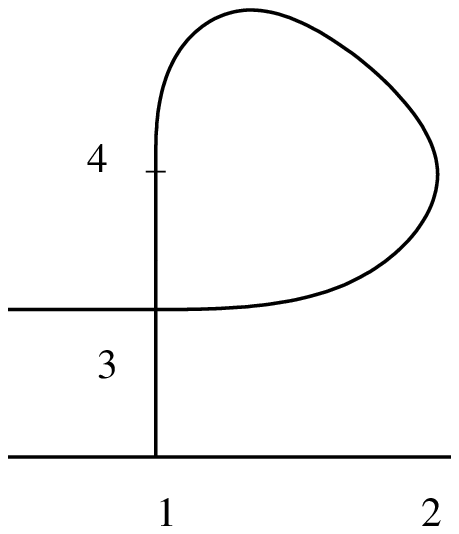}\hspace{1cm}\includegraphics[height=5cm]{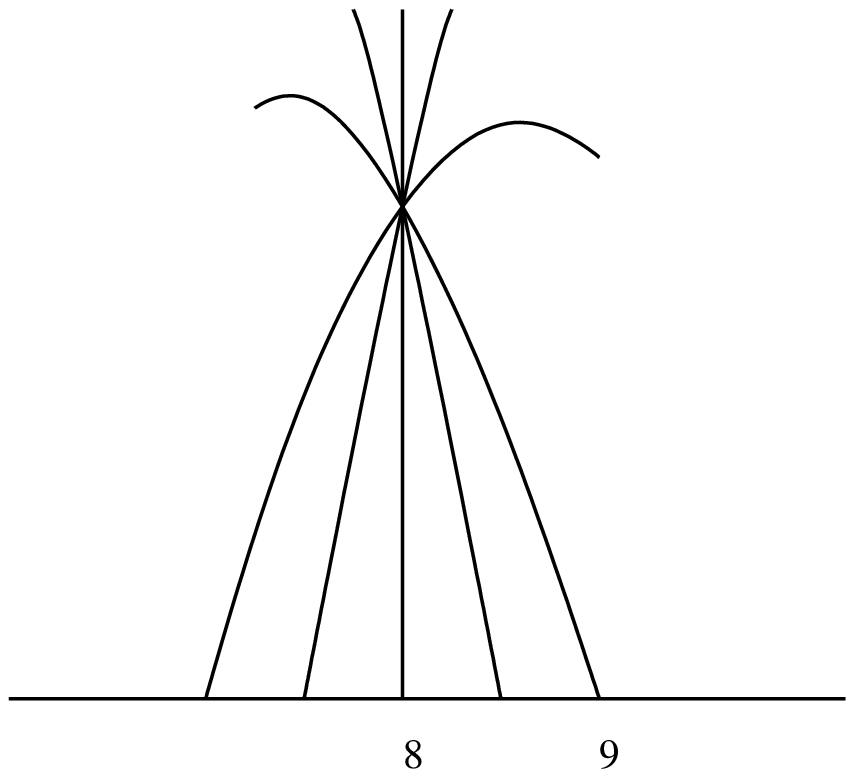}
\end{center}
\caption{Left: Self-intersection of a normal geodesic.
Right: Geodesics corresponding to focusing directions.}
\end{figure}

Later we will consider intersections of various geodesics on  $M$.
In these considerations we would like to avoid pathological cases that
may happen to long geodesics. The first case we  analyze is
a self-intersection of a geodesic.

\begin{lemma}\label{M 1} Let $\gamma_{z,\nu},\,z\in \p M$ be 
 the normal geodesic  and
\ba
\gamma_{z,\nu}(s_+)=\gamma_{z,\nu}(s_-),\quad s_+>s_-,
\ea
that is, $\gamma_{z,\nu}$ intersects itself.
Then $s_++s_->2\tau_R(z,\nu)$.
\end{lemma}

\noindent
{\bf Proof.} Assume that 
\beq
\label{28.4.11}
s_++s_-\leq 2\tau_R(z,\nu). 
\eeq
Then $s_-<\tau_R(z,\nu)$.
Let  $A=\gamma_{z,\nu}(s_-)$,
$B=\gamma_{z,\nu}(\tau_R(z,\nu))$ be points on $\gamma_{z,\nu}$, see Fig.\ 2, and denote by
$l_{BA}=s_+-\tau_R(z,\nu)$ the length of the "long" geodesic $\gamma_{z,\nu}
([\tau_R(z,\nu),s_+])$.
Then, using definition (\ref{28.4.2}) of $\tau_R$,
$s_- =\dist(z,A), \,  \tau_R(z,\nu) -s_- = \dist(A,B)$, so that 
the length of the broken geodesic $\gamma_{z,\nu}([0,s_+])\cup 
\gamma_{z,\nu}([0,s_-])$ from $z$ to $z$ is 
\ba
s_++s_-=\dist(z,A)+\dist(A,B)+l_{BA}+\dist(A,z).
\ea
Since 
$\gamma_{z,\nu}([s_-,\tau_R(z,\nu)])$ is the unique minimal geodesic
between its endpoints, 
$
l_{BA}> \dist(A,B) = \tau_R(z,\nu)-s_-.
$
Therefore,
\ba
s_++s_-
 >s_-+(\tau_R(z,\nu)-s_-)+(\tau_R(z,\nu)-s_-)+s_-
=2\tau_R(z,\nu),
\ea
which contradicts (\ref{28.4.11}).
\proofbox 

In the sequel, $\dist_{S}$ is the Sasakian distance on,
depending on the context,  $T\tilde{M}$ or $S\tilde{M}$,
see \cite{Sa}.

\begin{lemma}\label{M 2} Let $\e >0,\, z \in \p M$. There is 
$\delta = \delta(\e)>0$ such that
if
\bfo
(z_1,\xi_1)\, R_{2t}\,(z_2,\xi_2), \quad \hbox{i.e.} \,\,
\gamma_{z_1,\xi_1}(t_1)=\gamma_{z_2,\xi_2}(t_2), \,\,\,
t_1+t_2=2t,
\efo
with $t < \tau_R(z,\nu) +\delta$ and $\dist_{S}((z_i,\xi_i),\,(z,\nu)) < \delta$, $i=1,2$
then
\bfo 
|t-t_i| <\e,\quad i=1,2.
\efo
\end{lemma}

Note that the constant $\delta$  does not depend on $z\in \p M$.

\noindent
{\bf Proof.}
Assume the opposite, i.e., an existence of  points
$z^k \in \p M$,  $(z_i^k,\xi_i^k) \in \Omega_+,\, k=1,2,\, i=1,2,\dots$
and a parameter $\e >0$,
such that 
\bfo
\lim_{k\to \infty} \dist_{S}((z_i^k,\xi_i^k),\,(z^k,\nu^k))=0, 
\efo
\bfo
\gamma_{z_1^k,\xi_1^k}(t_1^k)=\gamma_{z_2^k,\xi_2^k}(t_2^k), \,\,\,
t_1^k+t_2^k=2t^k, \,\, \limsup_{k \to \infty} (t^k - \tau_R(z^k,\nu^k)) \leq 0,
\efo
with $t_1^k-t_2^k \geq 2\e$.
Using continuity arguments and  compactness of $\p M$ we have that 
there is a subsequence $k(p)$ with 
$z^{k(p)}\to z$, $t^{k(p)}_1\to t^+$,  $t^{k(p)}_2\to t^-$, and 
\bfo
\gamma_{z,\nu}(t_+) = \gamma_{z,\nu}(t_-), \quad t_+ + t_- \leq 2 \tau_R(z,\nu),
\,\,\,
t_+ - t_- \geq 2\e,
\efo
which contradicts Lemma \ref{M 1}.
\proofbox

Next we introduce  auxiliary functions $\mu_1(z)$, $\mu_2(z)$,
and $\tau_M(z)$, $z\in \p M$ with $\mu_1(z)$ and
$\mu_2(z)$ to be determined from the broken scattering relation.
{The function $\mu_1(z)$ tells when a normal geodesics sent from
 $z \in M$ exits $M$.  By the definition of the broken scattering relation, $R$,
 a point $(z, \xi) \in \Omega_+$ is in 
relation with itself, $(z, \xi)R_t(z, \xi)$, if and only if 
the geodesic $\gamma_{z, \xi}((0,t/2])$ on $\tilde M$ lies in $M^{\rm int}$.
This makes it possible to determine, for any $\gamma_{z, \xi},\, (z, \xi) \in \Omega_+$,
its arclength to the first hitting point to $\p M$. We denote this arclength 
by $\mu_1(z, \xi)$ and $\mu_1(z)=\mu_1(z, \nu)$.

The function $\mu_2(z)$ is an approximation
to $\tau_f(z)$. If we want to determine $\tau_f(z)$ we can argue as follows: 
assume that $s >\tau_f(z)$. }
Then the  normal geodesic
 $\gamma_{z, \nu}([0,s])$ 
is no longer a shortest path from $\gamma_{z, \nu}(s)$ to $\p M$
and there are
sequences $ z_n \to z,$  $z_n\not=z$, $ s_n \to \tau_f(z), t_n \to \tau_f(z)$ such that
\ba
%\label{12.5.0}
\gamma_{z, \nu}(s_n)=\gamma_{z_n, \nu_n}(t_n),\quad \nu_n=\nu(z_n).
\ea
In terms of the relation $R$, these imply that
\beq
\label{12.5.1}
& & \quad \quad(z,\nu)\,R_{T_n}\,(z_n, \nu_n), \quad T_n= t_n+s_n, \\ \nonumber
\noindent
&\hbox{with}& 
s_n \to \tau_f(z),\,t_n\to \tau_f(z), \, z_n \to z,\,\, \hbox{when} \,\, n \to \infty.
\eeq
Therefore, it  makes
 sense to try to find $\tau_f(z)$ using (\ref{12.5.1}). 
However, there are  two obstacles. 
First, it may happen that $\tau_f(z) \geq \mu_1(z)$.
Second,  having (\ref{12.5.1}) with $z_n \to z, \, T_n \to 2t$, we want to conclude
that $s_n \to t,\, t_n \to t$. To do so, we intend to use Lemma \ref{M 2}, which requires $t \leq \tau_R(z,\nu)$ which is not known. 
To avoid these difficulties, we will not determine
$\tau_f(z)$ but another function
$\mu_2(z)$ that is closely related to it.

\begin{definition}
\label{pseudofocal}
Consider the set $S(z)$ of those $s\in (0,\mu_1(z))$
 for which there are sequences $z_n \to z$, $z_n\in \p M$
$z_n \not= z$,  $T_n \to 2s$ such that
\beq
\label{12.5.3}
(z_n,\nu_n)\, R_{T_n}\,(z,\nu).
\eeq
Define $\mu_2(z) = \inf S(z)$, if $S(z)\not =\emptyset $ and 
 $\mu_2(z)=\mu_1(z)$ otherwise. 
\end{definition}

Observe that $\mu_2$ may be found from the broken scattering relation.

\begin{lemma}\label{lem: foc} Function $\mu_2:\p M\to \R_+$ 
satisfies
\beq\label{mu 2- cond}
\min(\mu_1(z),
\tau_f(z),\tau_R(z,\nu))
\leq \mu_2(z)\leq \min(\mu_1(z),\tau_f(z)).
\eeq
and $\tau_b(z)\leq \mu_2(z)$.
\end{lemma}

\noindent
{\bf Proof.}
The right inequality in (\ref{mu 2- cond}) follows from Definition
\ref{pseudofocal} and considerations before it.

To prove the left inequality of (\ref{mu 2- cond}), let us assume that there is 
$s<\min(\tau_f(z),\, \mu_1(z), \, \tau_R(z,\nu))$ which satisfies (\ref{12.5.3}). By Lemma
\ref{M 2}, applicable due to $T_n < 2 \tau_R(z,\nu)$ for large $n$,
we have 
\beq
\label{12.5.4}\quad\quad
\gamma_{z_n,\nu_n}(s_n)=\gamma_{z,\nu}(s_n'), \quad s_n\to s,\ \ \,s_n' \to s,
\quad z_n \to z,\ z_n\not= z.
\eeq
As $s < \tau_f(z)$, $\exp_{\p M}$ is a local diffeomorphism near $(z,s)$,
which contradicts (\ref{12.5.4}). This proves (\ref{mu 2- cond}).

Using definitions $\mu_1$ and $\tau_f$, we see by using
(\ref{24.5.7}) 
that 
\ba
\tau_b(z)\leq \min(\frac 12 \mu_1(z),
\tau_f(z),\tau_R(z,\nu(z))).
\ea
This yields $\tau_b(z)\leq \mu_2(z).$
\proofbox

Finally, we need a function $\tau_M(z)$ with
$\tau_M(z)>\tau_b(z)$
 having the property that, for $t<\tau_M(z)$ the 
geodesics sent back from a point 
$x=\gamma_{z,\nu}(t)$ hit the boundary $\p M$
near $z$ in a regular way.
Namely, we define
\ba
\tau_M(z)=\min\big(\mu_1(z),\tau_R(z,\nu(z))\big),\quad z\in \p M.
\ea
As $\tau_b(z)\leq \frac 12 \mu_1(z)$ we see by (\ref{24.5.7}) that 
$\tau_b(z)< \tau_M(z)$.

\subsection{Family of intersecting geodesics}

In this section we intend to use the broken scattering relation to
verify if a given family of  geodesics intersect
at one point.

Let $z_0\in \p M$,
$\nu_0=\nu(z_0)$, and $x_0 = \gamma_{z_0,\nu_0}(t_0),$ 
$ 0<t_0 <\tau_M(z_0)$. 
Denote  $\eta_0 = -
\gamma'_{z_0,\nu_0}(t_0)$.
Clearly, $\eta_0$ is the direction of the reverse geodesic,
$\gamma_{x_0,\eta_0}$ from $x_0$ to $z_0$.  By considering
Jacobi fields along this geodesic, we see that the exponential map,
$\exp_{x_0}: S_{x_0}\tilde{M} \times  \R_+ \to \tilde{M}$,  is a
local diffeomorphism near $(\eta_0,t_0)$.

{As $t_0 <\tau_R(x_0, \eta_0)$ and $\gamma_{x_0, \eta_0}(t_0)$ 
hits $\p M$ normally, all geodesics $\gamma_{x_0, \eta}$ hit $\p M$
transversally for $\eta \in S_{x_0}M$ close to $\eta_0$. They determine 
smooth functions $z(\eta),\, t(\eta)$ such that $\gamma_{x_0, \eta}(t(\eta))=z(\eta) \in \p M$.
Inverting these functions and using transversality, we obtain, in a neighborhood 
$U \subset \p M$ of $z_0$ a smooth section $\xi(z): U \to SU$ and a function $t(z)$
such that
}
\beq\label{eq: goal}
\gamma_{z,\xi(z)}(t(z))=x_0,\quad z\in U.
\eeq

In the following, our aim is to determine, 
using the broken scattering relation $R$,
whether, for a given 
 triple $\{U,\xi(\,\cdotp),t(\,\cdotp)\}$
of a neighborhood $U\subset \p M$ and functions $\xi(z)$ and $t(z)$,
there exists a point $x_0\in M$
such that $\gamma_{z,\xi(z)}(t(z))=x_0$ for all $z\in U$.

To this end, we notice that  property (\ref{eq: goal})  
implies
\beq
\label{12.5.10}
& &(z,\xi(z))\, R_{T(z)}\,(z_0,\nu_0), \quad (z,\xi(z))\, R_{T(z,z')}\,
(z',\xi(z')), \quad z,z' \in U, \\
\nonumber
& &T(z)=t(z)+t_0,\quad
T(z,z') = t(z)+t(z'),
\eeq 
for smooth $\xi(z),\, t(z)$. In addition,
 \beq
\label{12.5.11}
t(z_0) =t_0,\quad dt(z)|_{z_0}=0,\quad \xi(z_0)=\nu(z_0),  
\eeq
where the last properties follow from the fact that 
$\gamma_{x_0, \eta_0}$ is normal to $\p M$.
Here, $dt(z)=d_z t(z)$ is the differential of the function $t:U\to \R$.

These observations motivate the following definition:

\begin{definition}
\label{focusing} Let $z_0\in \p M$ and $t_0>0$.
Consider a family $ \mathcal{F}
(z_0,t_0)=\{U,\xi(\,\cdotp),t(\,\cdotp)\}$
where $U\subset \p M$ is a neighborhood of $z_0$,
$\xi:U\to SM$ is a smooth section, and $t:U\to \R$ is a smooth function.
We say that $\mathcal{F}(z_0,t_0)$ is a family of focusing directions if 
$\xi(z),\, t(z)$ satisfy  
conditions (\ref{12.5.10}) and (\ref{12.5.11}).
We then say that  the geodesics $\gamma_{z, \xi(z)},\ z\in U$ are the geodesics 
corresponding to family ${\mathcal F}(z_0,t_0)$.
\end{definition}

Note that the broken scattering relation $R$
determines if given $U$, $\xi(z)$, and $t(z)$
form a family of focusing directions. Our principal technical result in
this section shows  that 
the geodesics 
corresponding to a family of focusing directions
intersect at a single point.

\begin{theorem}\label{unique focusing} 
  Let $z_0 \in \p M,\, t_0 <\tau_M(z_0)$, and $\mathcal{F}
  (z_0,t_0)$ be a family of focusing directions. Then there is
a neighborhood  $\tilde{U} \subset U$ of $ z_0$
such that
 \ba
%\label{13.5.1}
\gamma_{z,\xi(z)}(t(z)) = \gamma_{z_0,\nu_0}(t_0), \quad
\hbox{for all} \,\, z \in \tilde{U}.  
\ea
\end{theorem}

\noindent
{\bf Proof.}  The proof of this result is rather long and will consist
of several steps and auxiliary lemmata.

\noindent {\it Step 1.}
We start with an observation that (\ref{12.5.10}) implies that,
 for any $z\in U$, there are $s(z), \hat{s}(z)\geq 0$ such that 
 \ba
%\label{12.5.12}
x(z)=\gamma_{z,\xi(z)}(s(z)) = \gamma_{z_0,\nu_0}(\hat{s}(z)), \quad
s(z)+\hat{s}(z)= T(z).  
\ea 
As $t_0 < \tau_R(z_0, \nu_0)$, by Lemma \ref{M 2} $s(z)\to t_0$, 
$\hat{s}(z) \to t_0$ when $z \to z_0$ and
\beq
\label{eq: 27 B}
s(z_0)=\hat s(z_0)=t_0.
\eeq
Next we show that $s(z), \hat{s}(z)$ are $C^{\infty}$-smooth near $z_0$ and 
\beq
\label{12.5.13}
ds(z)|_{z_0} = d\hat{s}(z)|_{z_0} =0.  
\eeq 
To this end, consider the
function $H(s,z)$,
 \ba
H(s,z) = \dist(\gamma_{z_0,\nu_0}(s), z) +s -T(z), \quad (s,z)
\in (t_0-\delta, t_0+\delta) \times U.
\ea 
As $t_0 <\tau_R(z_0,\nu_0)$, the function $H(s,z)$
is $C^\infty$-smooth a neighborhood of  $(t_0,z_0)$
and 
\ba
%\label{13.5.2}\quad\quad
H(t_0,z_0)=0,\quad \p_sH (t_0,z_0)= \p_s
\dist(\gamma_{z_0,\nu_0}(s),\,z_0)|_{t_0}+1=2.\hspace{-1cm}  
\ea
Making $U$ smaller if necessary,   the equation
$H(s,z)=0$
 has a unique solution $s=\tilde{s}(z)$ which is 
$C^{\infty}-$smooth in $U$ with $ \tilde{s}(z_0)=t_0$.
As also $s=\hat s(z)$ solves  $H(s,z)=0$,
we see  that 
 $\hat{s}(z)=\tilde{s}(z)$, $z \in U$.  It then
follows that $s(z)=T(z)-\hat s(z) \in C^{\infty}(U)$.

Let us differentiate the identity $H(\hat{s}(z),z)=0$ with respect to
$z$ at $z=z_0$. Due to (\ref{12.5.11}) and the fact that
$\gamma_{z_0,\nu_0}$ is normal to $\p M$, 
\bfo 
0 =d_zH(\hat{s}(z),z)|_{z_0} =
d_z \hat{s}\,|_{z_0} \cdot (\p_s \dist(\gamma_{z_0,\nu_0}(s),\,z_0)|_{s=t_0}+1)=2d_z \hat{s}\,|_{z_0} .  
\efo
Thus,
$d_z \hat{s}\,|_{z_0}=0$ and also
 $d_z s|_{z_0}
=d_z (T(z)-\hat s(z))|_{z_0}=0.$

\noindent {\it Step 2.} Consider the map $E \in C^{\infty}(U; SM)$,
\ba
E(z) = (x(z), \eta(z)):= \left(\gamma_{z,\xi(z)}(s(z)),\,
  -\gamma'_{z,\xi(z)}(s(z))\right), \quad E(z_0) =
(x_0,\eta_0).  \ea

\begin{lemma}\label{bijection}
  The map $dE|_{z_0}: T_{z_0} \p M \to T_{x_0,\eta_0}SM$ has the form
  \beq
\label{12.5.16}
dE|_{z_0} (v) =(0, \Theta v), \quad v \in T_{z_0} \p M , \eeq where we
identify $T_{x_0,\eta_0}SM \approx T_{x_0}M \times
T_{\eta_0}(S_{x_0}M)$.  Furthermore, $\Theta: T_{z_0} \p M \to
T_{\eta_0}(S_{x_0}M)$ is bijective.
\end{lemma}

\noindent{\bf Proof of Lemma \ref{bijection}.} 
As $x(z) = \gamma_{z_0,\nu_0}(\hat{s}(z))$, it
follows from (\ref{12.5.13}) that $dx|_{z_0} =0$, i.e.,  $dE|_{z_0}$ is
of  form (\ref{12.5.16}). To show that $\Theta$ is bijective, observe
that \beq
\label{12.5.17}
 \exp_{x(z)}(s(z)\eta(z)) =
z, \quad z \in U.  \eeq
Let us denote $\hbox{Exp}(x,\xi)=\exp_x \xi,\, (x, \xi) \in T{\tilde M}$. By differentiating 
both sides of (\ref{12.5.17}) with respect to $z$ 
and using
$dx|_{z_0} =0$, we obtain
\ba
%\label{12.5.18}
d_\xi\hbox{Exp}|_{(x_0,t_0 \eta_0)}\bigg(s(z_0)\Theta \zeta+ 
(ds|_{z_0}\zeta)\eta(z_0)\bigg)=\zeta
\ea 
for any $\zeta\in T_{z_0} \p M.$ Using that $s(z_0)=t_0, ds|_{z_0}=0$, we get 
\bfo
d_{\xi}\exp_{x_0}|_{\xi=t_0\eta_0}(t_0 \Theta \zeta)=\zeta,
\efo
which implies that $\Theta: T_{z_0}\p M
\to T_{\eta_0}(S_{x_0}M)$ is bijective.  \proofbox

\noindent {\it Step 3.}
Our further considerations are based on the analysis of 
the
intersection of a single geodesic and 
the geodesics corresponding to a family of focusing directions.

\begin{lemma}\label{non-inters} 
Let
$z_0\in \p M$ and   $\mathcal{F} (z_0,t_0)=\{U,\xi(\,\cdotp),t(\,\cdotp)\}$,
$t_0 <\tau_M(z_0)$
be a  family  of 
focusing directions. Let $\gamma(\tau)$ be another geodesic in 
${M}$ which intersects $\gamma_{z_0, \nu_0}$,
\beq
\label{16.5.0}\quad\quad
\gamma(0)=\gamma_{z_0,\nu_0}(r_0),\quad
\gamma'(0) \neq \pm \gamma'_{z_0,\nu_0}(r_0),\quad r_0<\tau_M(z_0).
  \eeq 
Assume, in addition, that all geodesics $\gamma_{z,\xi(z)}$ corresponding to
$\mathcal{F} (z_0,t_0)$
intersect $\gamma$
near $y_0$, i.e., 
\beq
\label{20.5.1}
\gamma_{z,\xi(z)}(r(z)) = \gamma(\tau(z)), 
\eeq 
where
 $0<r(z) \leq r_1<\tau_M(z_0)$ and
$|\tau(z)| \leq i_1<\hbox{inj}\,(M)$.
Then $r_0=t_0$.
\end{lemma}

\noindent
{\bf Proof of Lemma \ref{non-inters}.}
Denote $y_0=\gamma_{z_0,\nu_0}(r_0)$.
First we show that $r(z)$ is continuous
at $z_0$. If this is not true, there would be another intersection of
$\gamma_{z_0,\nu_0}$ and $\gamma$, 
\bfo 
\gamma_{z_0,\nu_0}(r')
=\gamma(\tau'), \quad r' \leq r_1,\
r' \not=r_0,\ |\tau'| < \hbox{inj}\,(M).
\efo 
This  leads to a contradiction as both $\gamma([0,\tau'])$
and $\gamma_{z_0,\nu_0}([r_0,r'])$ are unique minimal geodesics
between their endpoints. 
Thus $r(z)$ is continuous at $z_0$.

To prove the claim, we assume  that $r_0 \neq t_0$. 
Our next goal is to show that the map
$\Psi: U \times \R_+\to M$, 
\bfo 
\Psi(z,r) =\exp_z(r\xi(z)) 
\efo 
is a local
diffeomorphism near $(z_0,r_0)$, see the right part of Fig.\ 3.
{Indeed, as $t_0, r_0 < \tau_R(z_0, \nu_0)$, the map
$\exp_{x_0}$ is a local diffeomorphism near $(t_0-r_0)\eta_0$,
where $x_0=\gamma_{z_0, \nu_0}(t_0),\, \eta_0 =- \gamma'_{z_0, \nu_0}(t_0)$.
Thus,
}
 \ba
%\label{16.5.1} 
d \exp_{x_0}|_{(t_0-r_0)\eta_0}: T_{(t_0-r_0)\eta_0}(T_{x_0}M)\to
T_{y_0}M 
\ea 
is bijective.  Using the definitions for $s(z), \, E(z)=(x(z),\eta(z))$ introduced
earlier we have 
\ba
%\label{16.5.2}
\Psi(z,r) = \gamma_{E(z)}(s(z)-r) = \exp_{x(z)}((s(z)-r)\eta(z)).
\ea
By (\ref{eq: 27 B}) and (\ref{12.5.13}), $ds(z)|_{z_0}=0$ and $s(z_0)=t_0$,
which
together with
 (\ref{12.5.16}) imply that 
\ba
%\label{16.5.3}
d\Psi|_{(z_0,r_0)}(\zeta,\rho)
= d \exp_{x_0}|_{(t_0-r_0)\eta_0}
((t_0-r_0)\Theta \zeta-\rho\eta_0)
\ea
for $\zeta\in T_{z_0}\p M$ and $\rho \in \R$.
Thus, by Lemma \ref{bijection} and bijectivity of 
$d
\exp_{x_0}|_{(t_0-r_0)\eta_0}$, 
\bfo 
d\Psi|_{(z_0,r_0)}:T_{z_0}\p M \times \R \to T_{y_0}M 
\efo 
is bijective, i.e., $\Psi$ is a local diffeomorphism
near $(z_0,r_0)$. 
%{\newtextt Thus, making $U$ smaller if
%necessary, we can assume that the map $\tilde{\Psi}(z) =\Psi(z,r(z))$
%defines a  diffeomorphism 
%$\tilde \Psi:U\to \tilde \Psi(U) \subset \Sigma$.}

Now, let $\Sigma$ be an $(n-1)-$dimensional submanifold which contains 
a part $\gamma(-\e, \e)$ of $\gamma$ near
$y_0$ and is transversal to $\gamma_{z_0,\nu_0}$ at $y_0$, see Fig.\ 3, 
the existence of such submanifold guaranteed by   (\ref{16.5.0}).
Introducing the boundary normal coordinates
$(w,n)$ associated to $\Sigma$, with $n=0$ on $\Sigma$, we rewrite $\Psi$ 
in these coordinates as 
\bfo
\Psi(z,r) = (w(z,r),\,n(z,r)).  
\efo 
By transversality, 
$
\frac{\p  n}{\p r}(z_0,r_0) \neq 0.
$ 
{This implies 
that
for any $z$ near $z_0$ the equation 
$n(z,r)=0$ for has a unique solution 
$r={\hat r}(z)$.}
Moreover,  ${\hat r}(z_0)=r_0$ and the function ${\hat r}(z)$ is smooth in 
a neighborhood of $z_0$.

\begin{figure}[htbp]\label{Figure of center x_0}
\begin{center}
\psfrag{1}{$\p M$}
\psfrag{2}{$z_0$}
\psfrag{3}{$\gamma$}
\psfrag{4}{$\Sigma$}
\psfrag{5}{$x_0$}
\psfrag{6}{$p$}
\includegraphics[width=6cm]{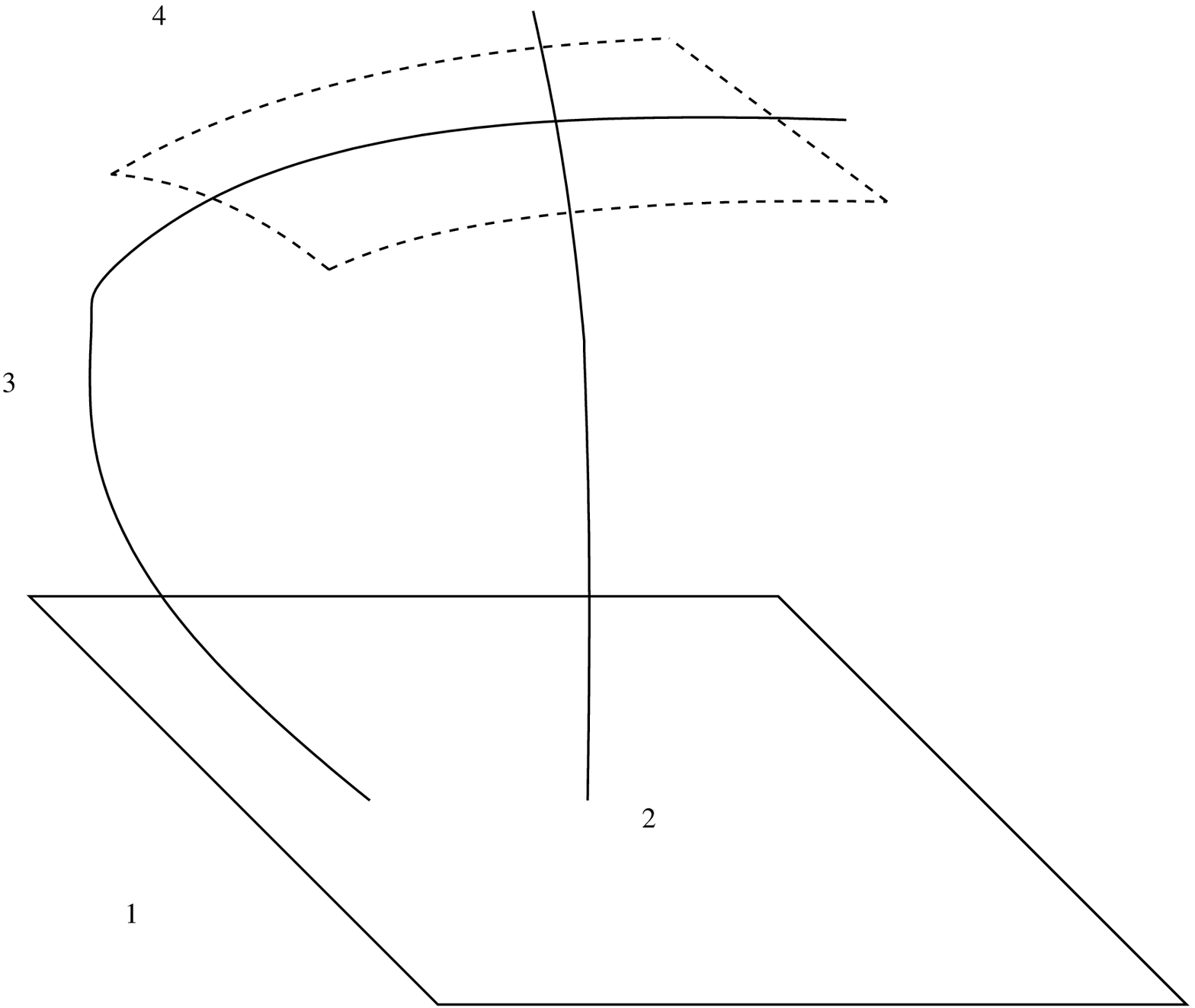}
\includegraphics[width=6cm]{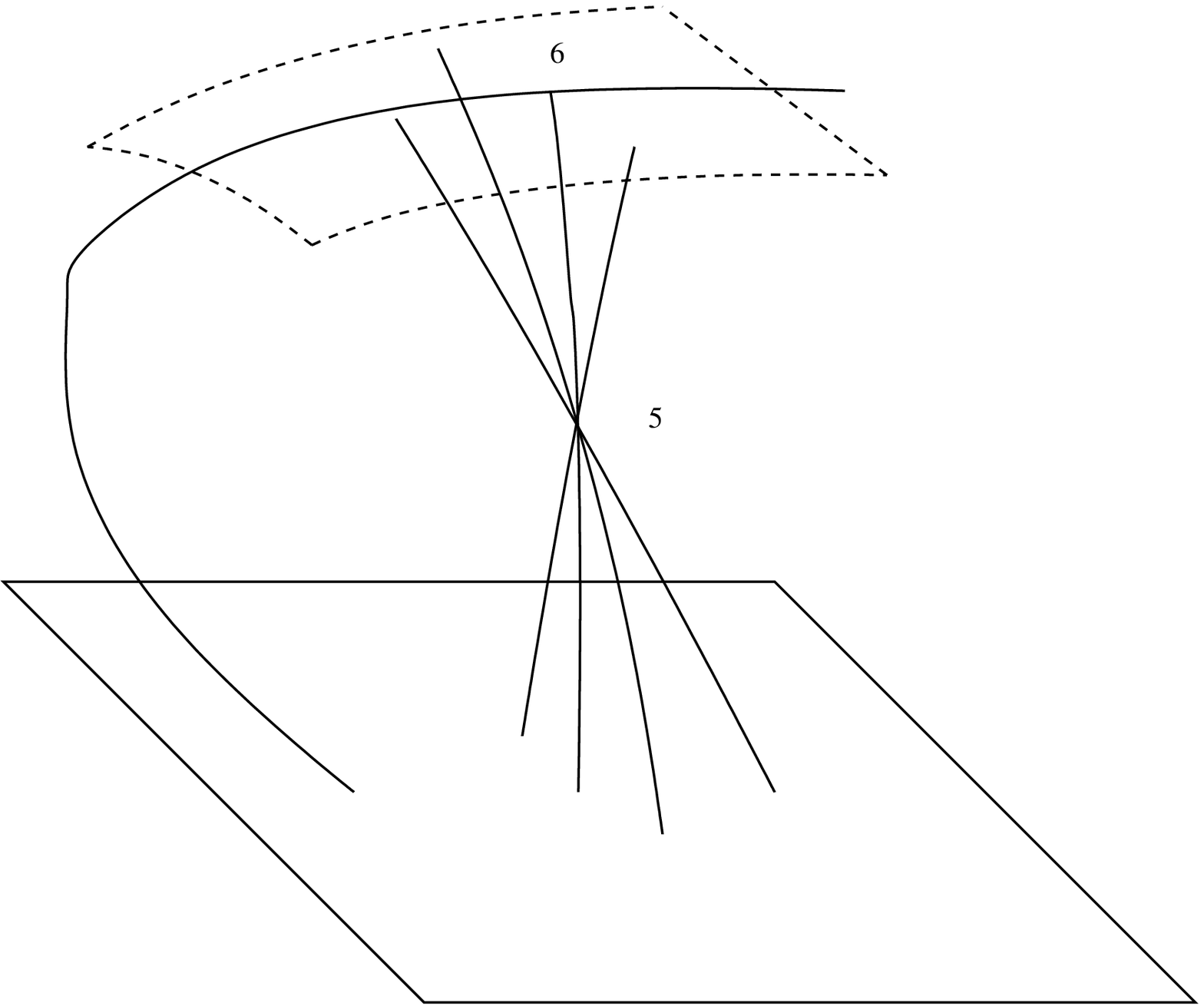}
\end{center}
\caption{Left: Submanifold $\Sigma$ contains geodesic $\gamma$
  and is transversal to $\gamma_{z_0,\nu}$. Right: Geodesics
corresponding to $\F(z_0,t_0)$ almost intersect  at the point 
$x_0=\gamma_{z_0,\nu}(t_0)$ and define coordinates near 
$p=\gamma_{z_0,\nu}(r_0)$.}
\end{figure}

Now  $r(z)$ and $\hat r(z)$ are continuous at $z_0$ and 
they both solve the equation $n(z,r)=0$.  Thus, there is
a neighborhood $\tilde U\subset U$ of $z_0$ such that 
$\hat{r}(z)=r(z)$ for $z\in \tilde U$. 
As also $\Psi$ is a local diffeomorphism, we see that if $\tilde U$ is small enough, then
$\tilde \Psi:\tilde U\to \tilde \Psi(\tilde U) \subset \Sigma$, where $\tilde{\Psi}(z) =\Psi(z,r(z))$,
is a diffeomorphism of $(n-1)$-dimensional submanifolds.
%Above we saw that $\tilde \Psi:U\to \tilde \Psi(U) \subset \Sigma$
%is a diffeomorphism.
On the other hand, condition (\ref{20.5.1}) 
implies that $\tilde \Psi(\tilde U)\subset  \gamma(-\e, \e)$.
 As $\gamma(-\e, \e)$ is a one-dimensional submanifold
of $\Sigma$, we get a contradiction for $n \geq 3$.
Thus, $r_0=t_0$.
\proofbox

\noindent {\it Step 4.}
Let
$0<\e<\frac 14\min (\hbox{inj}\,(M),\tau_R(z_0,\nu)-t_0)$
and $0<\delta<\delta(\e)$ where  $\delta(\e)$ is defined in Lemma \ref{M 2}. 
 We choose a neighborhood $\tilde U\subset U$ of $z_0$ so that 
\ba
|t(z)-t_0|<\e\quad\hbox{and}\quad 
d_{S}((z,\xi(z)),(z_0,\nu_0))<\delta\quad\hbox {for }z\in \tilde U.
\ea

By Definition \ref{focusing},
 there exist functions $s_1(z,z'),\, s_2(z',z)>0,\,z,z' \in \tilde U,$ 
such that
 \ba
\gamma_{z,\xi(z)}(s_1(z,z')) = \gamma_{z',\xi(z')}(s_2(z',z)), \quad s_1(z,z')
+s_2(z',z) =t(z)+t(z').  
\ea
By Lemma \ref{M 2}, these imply that
\beq
\label{25.7.13}
|t_0-s_1(z,z')| <2\e,\quad |t_0-s_2(z',z)| <2\e.
\eeq
Consider a geodesic $\gamma(s)= \gamma_{z',\xi(z')}(s+s_2(z',z_0))$
 for some
 fixed $z' \in \tilde U$, $z'\not =z_0$. It
follows from (\ref{25.7.13}) that Lemma \ref{non-inters}
 is applicable to the family
${\mathcal F}(z_0,t_0)$ and the geodesic 
$\gamma$ with $r_1 = \tau_R(z_0,\nu_0)-2\e,\,
i_1=2\e$.  Thus,
$\gamma_{z',\xi(z')}$ and $ \gamma_{z_0,\nu_0}$
intersect at $x_0=\gamma_{z_0, \nu_0}(t_0)$.
As $z' \in \tilde U\setminus \{z_0\}$ is arbitrary, 
all geodesics corresponding
to family ${\mathcal F}(z_0,t_0)$ with a starting
point $z'\in \tilde U$ intersect in 
$x_0.$ 
\proofbox

Later on we will need  the following modification of Lemma
\ref{non-inters} which do not require that all geodesics of ${\mathcal
  F}(z_0,t_0)$ intersect $\gamma$ near $y_0$.

\begin{lemma}\label{non-inters, hausdorff}
Let
$z_0\in \p M$ and   $\mathcal{F} (z_0,t_0)=\{U,\xi(\,\cdotp),t(\,\cdotp)\}$,
$t_0 <\tau_M(z_0)$
be a  family  of 
focusing directions.
 Let $\gamma(\tau)$ be another geodesic in 
${M}$ which intersects all geodesics $\gamma_{z,\xi(z)}$ corresponding to
$\mathcal{F} (z_0,t_0)$,
\ba
\gamma_{z,\xi(z)}(r(z)) = \gamma(\tau(z)), 
\ea
where
 $0<r(z) \leq r_1<\tau_M(z_0)$ and
$|\tau(z)| \leq L$, where $L>0$ is arbitrary.
%\HOX{I've added details to the formulation, please look at}
Assume, in addition, that $h(z)=r(z)+\tau(z)$ is continuous.
Then $\gamma_{z,\xi(z)}(t(z))=\gamma(h(z_0)-t_0)$
when $z$ is sufficiently close to
$z_0$, i.e.,
all geodesics intersect at the same point.

%Let conditions of Lemma \ref{non-inters} be valid {with
%the condition $|\tau(z)| <i_1$
%replaced by $|\tau(z)| \leq L$, where $L>0$ is arbitrary, 
%{\newtextt and assuming
%that $h(z)=\tau(z)+r(z)$ is continuous.
%Then $\gamma_{z_0,\xi_0}(t_0)=\gamma(h(z_0)-t_0)$.}
%}
\end{lemma}

\noindent 
{\bf Proof.}  We first show that there is only a finite number of
intersections of $\gamma_{z_0,\nu_0}((0,r_1))$ with $\gamma([-L,L])$.
Let $\tau_1, \dots, \tau_N \in [-L,L]$ and $r^1_0,\dots, r^N_0 \in (0,r_1)$
define the points of the intersection, 
\bfo
\gamma_{z_0,\nu(z_0)}(r^j_0)= \gamma(\tau_j).  
\efo 
As all  geodesics in balls of radius $\hbox{inj}\,(M)$ are shortest and
$r^j_0 \leq r_1$ with
 $\gamma_{z_0,\nu_0}([0,r_1])$ being the shortest between its endpoints, 
 \bfo 
N \leq
\left[ \frac{2L}{\hbox{inj}\,(M)} \right] +1,  
\efo
where $[t]$ denotes the  integer part of $t\in \R$.

Let $0<\e<\frac 12\hbox{inj}\,(M)$ and $U(\rho)=\p M \cap 
{B(z_0,\rho)}$,
where $B(z_0,\rho)\subset M$ is the ball with center $z_0$ and radius $\rho$.
Then there is $\rho_0>0$ such that 
\bfo 
\min _{1\leq j \leq N} |r(z) -
r^j_0| < \e,\quad \hbox{for }z\in U(\rho_0). 
 \efo 
Indeed, otherwise there is a sequence $z_n \to
z_0$ 
with $r(z_n) \to \tilde{r} < \tau_R(z_0,\nu(z_0))$ and $
\tau(z_n) \to \tilde{\tau}, \, |\tilde{\tau}| \leq L,$ such that 
\bfo
\gamma_{z_0,\nu_0}(\tilde{r})= \gamma(\tilde{\tau}), \quad
\tilde{r} \neq r^j_0, \,\, j=1,\dots, N, 
\efo 
which is a contradiction.

{\newtextt For $0<\rho<\rho_0$, denote
\ba 
V_j(\rho) = \{z \in U(\rho): \, \gamma_{z,\xi(z)}(r)
=\gamma(\tau),\ r(z)+\tau(z)=h(z),\ |r - r^j_0| \leq \e \}.
\ea
 Sets $V_j(\rho)$ are relatively closed  $U(\rho)$ and, therefore,
measurable on $\p M$. 
As $\bigcup_{j=1}^N { V}_j(\rho) = U(\rho)$,
we see that for some $j$ the set $V_j(\rho)$ has
non-zero $(n-1)$-dimensional measure.
 However,
if $r_0^j\not =t_0$, 
the same
considerations as in the proof of Lemma \ref{non-inters}, by
replacing $r_0$ by $r_0^j$ and using a relatively open neighborhood
 $\tilde U\subset V_j(\rho)$ of $z_0$,
show that the set ${ V}_j(\rho)$ has $(n-1)-$dimensional measure 
equal to $0$ when   $\rho>0$ is small enough.
This  shows that there are $j$ and $\rho>0$ such that 
$r_0^j=t_0$ and $U(\rho)\setminus V_j(\rho)$
has  $(n-1)-$dimensional measure 
equal to $0$. Thus $ V_j(\rho)$ is dense in  $U(\rho)$. As $\e>0$ is arbitrary, the
 continuity of the geodesic flow shows that
$\gamma_{z_0,\xi_0}(t_0)=\gamma(h(z_0)-t_0)$.
}
Together with Theorem \ref{unique focusing}  this completes the proof.
\proofbox 

In the following we say that two geodesics $\mu(t)$ and $\tilde \mu(t)$
 coincide if 
$\mu(t_1)=\tilde \mu (t_2)$ and
$\mu'(t_1)= \pm \tilde \mu '(t_2)$ for some $t_1,t_2\in \R$.
Note that this is equivalent to 
$\mu(t)=\tilde \mu(a+t)$ or $\mu(t)=\tilde \mu(a-t)$ for 
all $t$ in a non-empty open interval and $a\in \R$.

\subsection{Reconstruction of the boundary cut locus distance}

\begin{lemma}\label{lem: bou} The boundary, $\p M$, and 
the broken scattering relation, $R$, 
 determine the boundary cut locus distance
$\tau_b(z)$, $z\in \p M$.
\end{lemma}

\noindent
{\bf Proof.} We recall that for $t_0 < \tau_b(z_0)$ the point $z_0$ in the 
unique point of $\p M$ closest  to $x_0=\gamma_{z_0,\nu_0}(t_0)$. 
On the contrary, when
$t_0 > \tau_b(z_0)$ there is another point $w \in \p M$ with
$\dist(\gamma_{z_0,\nu_0}(t_0),\, w) < t_0.$
What is more, 
considerations in the beginning of Section 2.2 show the existence of a family 
${\mathcal F}(z_0,t_0)$ of focusing directions
for $t_0 < \tau_M(z_0)$. 
Recall that
$\tau_b(z_0)<\tau_M(z_0)$.

Thus, when $\tau_b(z_0) < t_0< 
\tau_M(z_0)$,
there is a family ${\mathcal F}(z_0,t_0)=\{U,\xi(\cdotp),t(\cdotp)\}$
of focusing directions,
 a point $w \in \p M, \, w \neq z_0$, and
$s_0 <t_0$ such that
\beq
\label{24.5.6}
(z,\xi(z)) \, R_{t(z)+s_0} \, (w, \nu(w)), \quad z \in U.
\eeq

Our next aim is to show that when $t_0 < \tau_b(z_0)$, there are 
no $w \in \p M$ and ${\mathcal F}(z_0,t_0)$
satisfying (\ref{24.5.6}) with $s_0 <t_0$.

Assuming the opposite, there is a neighborhood $U \subset \p M$ 
of $z_0 $ and a function
$r(z)$ with
\beq
\label{24.5.3}
\gamma_{z,\xi(z) }(r(z)) = \gamma_{w,\nu(w)}(t(z)-r(z)+s_0),\quad z\in U.
\eeq
Next we prove that
\beq\label{eq: A1}
r_0=\limsup_{z \to z_0} r(z) \leq t_0.
\eeq
Assume that  (\ref{eq: A1}) is not  true. Then there is a sequence $z_n \to z_0$ with $r(z_n) \to r_0 >t_0$. By the
continuity of the exponential map, it follows from (\ref{24.5.3}) that
$
\gamma_{z_0,\nu_0}(r_0) = \gamma_{w,\nu(w)}(t_0-r_0+s_0).
$
Thus, by the triangle inequality,
\ba
%\label{24.5.5}
& &\dist\left(w,\,\gamma_{z_0,\nu_0}(t_0)\right)\\
\nonumber& &\leq
 \dist(w, \gamma_{w,\nu(w)}(t_0-r_0+s_0))
+\dist(\gamma_{z_0,\nu_0}(r_0),\gamma_{z_0,\nu_0}(t_0))
\\
\nonumber & &
\leq (t_0-r_0+s_0)+(r_0-t_0)\leq s_0 <t_0,
\ea
which contradicts the definition (\ref{28.4.8}) of $\tau_b$.
Thus (\ref{eq: A1}) is valid.

{Therefore, by making $U$ smaller if necessary, we have
\bfo
r(z) <\tau_M(z_0), \quad z \in U.
\efo
 Assume first that  geodesics 
$\gamma_{z_0,\nu_0}$ and $ \gamma_{w,\nu(w)}$
do not coincide. 
Applying
Lemma \ref{non-inters, hausdorff} with 
$\gamma(\tau) =  \gamma_{w,\nu(w)}(t_0+s_0-r_0+\tau)$
 and $L=2t_0$, we obtain
$
\gamma_{z_0,\nu_0}(t_0)= \gamma_{w,\nu(w)}(s_0).
$
As $s_0 <t_0$ this  contradicts with
 the definition of $\tau_b$.}
If  
$\gamma_{z_0,\nu_0}$ and $\gamma_{w,\nu(w)}$ coincide, condition
$w\not=z_0$ implies that
$\gamma_{z_0,\nu(z_0)}(t_0+s_0)=w$.
Then we would have 
$\dist(x_0,\p M)\leq \dist(x_0,w)\leq s_0<\tau_b(z_0)$, that is
not possible.

Finally, by Lemma \ref{lem: foc} the relation $R$ determines the function $\mu_2(z)$
satisfying  $\tau_b(z)\leq \mu_2(z)$.
Let $J(z_0)$ be the set of those  $t_0\in [0,\mu_2(z_0)]$
for which there are $w \in \p M$,  $s_0 <t_0$, and ${\mathcal F}(z_0,t_0)$
satisfying (\ref{24.5.6}). If $\tau_b(z_0)<\mu_2(z_0)$,
we see that $(\tau_b(z_0),\mu_2(z_0))\subset J(z_0)$.
Thus we can determine $\tau_b(z_0)$ by setting $\tau_b(z_0)=\inf J(z_0)$ if $J(z_0)\not =\emptyset$
and $\tau_b(z_0)=\mu_2(z_0)$  otherwise.
\proofbox 

\subsection{Boundary distance representation of $(M,g)$.}

Next we  construct
of isometry type of manifold $(M,g)$ by showing that the broken scattering
 relation, $R$, determines 
the boundary distance representation 
${\mathcal R}(M)$ of $(M,g)$ that is the set
\ba
%\label{6.6.1}
{\mathcal R}(M) =\{r_x:\ x \in M\} \subset C(\p M),
\ea
where $r_x:\p M\to \R$ are the boundary distance functions
\ba
r_x(z) = \dist(x,z), \quad z \in \p M.
\ea
It is well-known, e.g.  \cite{AK2LT,KaKuLa,Ku} 
that the set ${\mathcal R}(M)$ possesses a 
natural structure of a Riemannian manifold with the map
\bfo
{\mathcal R}:M \to {\mathcal R}(M), \quad {\mathcal R}(x)=r_x(\cdot),
\efo
being an isomorphism. What is more, this metric structure can be identified just from the 
knowledge of the set ${\mathcal R}(M)$. An additional advantage of dealing with 
${\mathcal R}(M)$ is the existence of a stable procedure to construct a metric approximation,
in the Gromov-Hausdorff topology,  to
$(M,g)$ given an approximation to ${\mathcal R}(M)$ in the Hausdorff topology 
on $L^{\infty}(\p M)$,  \cite{K2L}.
To construct ${\mathcal R}(M)$, we assume that the function $\tau_b$ is already known.
 We start with finding $\dist_{\p M}$ on $\p M$ which is inherited from 
$(M, g)$. We define that 
$\dist_{\p M}(z_1,z_2)=\infty$ when $z_1$ and $z_2$ lie on different components of $\p M$.

\begin{lemma}  The boundary, $\p M$, and the broken scattering
 relation, $R$, determine,
 for any $z_1, z_2\in \p M$, 
 the distance  $\dist_{\p M}(z_1,z_2)$
along   $\p M$.
\end{lemma}
 
\noindent
{\bf Proof.} It is enough to consider the case when
$z_1$ and $z_2$ are in the same  component of  $\p M$.

Using  
boundary normal coordinates, we see that there is $\e_0>0$ and $c_0 >0$
 such that
\beq\label{eq: abs}
|\dist(y_1, y_2)-\dist_{\p M}(y_1, y_2)| \leq c_0 \e^{3/2},
\eeq
if $\dist_{\p M}(y_1, y_2)\leq \e ^{3/4}, \, \e < \e_0$.
{\newtextt Let $x_2=\gamma_{y_2, \nu_2}(\e^{5/4})$. Making $\e_0>0$ smaller
if necessary, we see that}
there is a unique shortest geodesic in $M$,
$\gamma_{y_1, \xi_1}$, with $(y_1, \xi_1) \in \Omega_+$, from $y_1$ to $x_2$.
Moreover, using again  
boundary normal coordinates, we see that
\beq\label{eq: abs2}
|\dist(y_1, x_2)+\dist(x_2, y_2)-\dist_{\p M}(y_1, y_2)|\leq c_1 \e^{5/4}.
\eeq
Let $\mu=\mu([0,l])$ be a shortest geodesic of $\p M$ from $z_1$
to $z_2$. Let $N\in \Z_+$,
$\e = l/N$ and $y_j=\mu(\e j)$, 
$j=0,\dots,N$.  Define $x_j= \gamma_{y_j, \nu_j}(\e^{5/4})$
and associate with each $j=1,\dots, N$ a broken geodesic $\alpha_j$ 
which is the union of the geodesic from $y_{j-1}$ to $x_j$ and from $x_j$ to $y_j$.
Inequality (\ref{eq: abs2}) implies that 
{\newtextt if $N\to \infty$, then}
\beq\label{eq: final1}\quad\quad
|\dist_{\p M}(z_1, z_2)- \sum_{j=1}^N \left( \dist(y_{j-1}, x_j)+\dist(y_j, x_j)  \right)|
\leq c_2 \e^{1/4} \rightarrow 0,\hspace{-1cm}
\eeq

Motivated by this, define for $N\in \Z_+$ and $\e=1/N$
\ba
d_{N}(z_1, z_2)=\inf\,
\sum_{j=1}^{N} s_j,
\ea
where the infimum is taken over the points
$y_j\in \p M$, $j=0,1,\dots,N, \, y_0=z_1, y_N=z_2$,
which satisfy the following condition:
{\newtextt
For any $j=0, \dots, N-1$,
there are  $\eta_{j} \in S_{y_{j}}M,\, (\nu_{j}, \eta_{j})_g >0$
and positive $s_j <\e^{3/4}$ such that
\ba
\bigg( (y_{j},\eta_j),(y_{j+1},\nu(y_{j+1})), s_j\bigg)\in  R,
\quad j=0,1,\dots, N-1.
\ea
Using (\ref{eq: abs}) we see that 
$d_N(z_1,z_2)\geq \dist_{\p M}(z_1,z_2)-c_3\e^{1/2}$. On the 
other hand, as we saw in (\ref{eq: final1}),
 there are $y_j,\eta_j$, and $s_j$ such that }
%$s_j<\dist_{\p M}(y_{j-1}, y_j)+c \e^{5/4}$. 
%This implies
\bfo
|\dist_{\p M}(z_1, z_2)-d_N(z_1, z_2)| \leq c_4 \e^{1/4}= c N^{-1/4} \rightarrow 0,
\quad \hbox{when}\ \  N \to \infty.
\efo
Thus we get that
\ba
\dist_{\p M}(z_1,z_2)=\lim_{N\to \infty} d_{N}(z_1, z_2).
\ea
\proofbox

Next we determine the distance between boundary points with respect to
the metric $g$ in $M$.

\begin{lemma}  
\label{boundary distance} The boundary, $\p M$, and the
 broken scattering
 relation, $R$, determine the distance function $\dist(x_1,x_2)$
for $x_1,x_2\in \p M$
\end{lemma}

\noindent
{\bf Proof.} 
{\newtextt
By \cite{Al}, for any $x_1,x_2 \in \p M$ a shortest path connecting them
is a  $C^1-$path. Let $x(s)$,
$s \in [0,l]$, $l=\dist(x_1,x_2)$, $x(0)=x_1,$ $ x(l)=x_2$
be such a shortest path,  parameterized by the arclength, that connects
$x_1$ to $x_2$ in $M$.
Moreover, by \cite{Al} it holds that
 if $x(s) \in M^{\rm int}$ for $ s \in (a,b)$, then $x((a,b))$ is a shortest
geodesic between $x(a)$ and $x(b)$ in $M$.

Clearly, the set of $s \in [0, l]$ such that
$x(s) \in M^{\rm int}$ is open. By (\ref{eq: abs}), for any $\e >0$ there is a 
finite number points $a_i,$ $i=1, \dots, p$, $a_{p+1}=l$,
and $ b_i,$ $i=1, \dots, p$ with
$0\leq a_1 <b_1 \leq a_2 \dots <b_p\leq a_{p+1}=l$
such that $z_i=x(a_i),y_i=x(b_i)\in \p M$ and 
\beq\label{eq: sum 1}
\dist(x_1, x_2)
\leq \dist_{\p M}(x_1, z_1)+ 
%\dist_{\p M}(z_{p+1},x_2)+  
\quad \quad \quad\quad \quad \quad\quad \quad   \\
\nonumber
\quad\quad+
\left(\sum_{i=1}^{p} 
\dist(z_i, y_i)+
\dist_{\p M}(y_i, z_{i+1})\right)
\leq \dist(x_1, x_2)+\e\hspace{-2cm}
\eeq
and there are shortest paths $\gamma_{z_i,\eta_i}([0,l_i])$ in $M$ 
of length $l_i=b_i-a_i$
from $z_i$ to $y_j$ that satisfy 
 $\gamma_{z_i,\eta_i}((0,b_i-a_i))\subset M^{\rm int}.$ 
Next we will relate (\ref{eq: sum 1}) to the broken geodesic
relation. Recall that relation $R$ involved
broken geodesics that start and end non-tangentially to the boundary. 
Because of this, we consider for tangential $\eta_i$ the vector
$\xi_i=(1-h)^{1/2}\eta_i+h^{1/2}\nu(z_i)\in S_{z_j}M$.
If  $\eta_i$ is  non-tangential, we set $\xi_i=\eta_i$.
When $h>0$ is small enough and $s_i<l_i$ is sufficiently
close to $l_i$, we have that 
 $\gamma_{z_i,\xi_i}((0,s_i])\subset M^{\rm int}$,
and the closest boundary point to $\gamma_{z_i,\xi_i}(s_i)$,
denoted  $\tilde y_i$,  satisfies
\ba
\dist(\gamma_{z_i,\xi_i}(s_i),\tilde y_i)<\frac \e p,\quad
\dist_{\p M}(\tilde y_i,y_i)<\frac \e p.
\ea
Consider the broken geodesic from $z_i$ to $\tilde y_j$ 
which is the union of the geodesic from $z_i$
 to $\gamma_{z_i,\xi_i}(s_i)$ and from $\gamma_{z_i,\xi_i}(s_i)$ to 
$\tilde y_j$. It  has the length 
$t_i\leq l_i+\e/p$ and non-tangential starting and ending directions.
Thus $(z_i,\xi_i)R_{t_i}(\tilde y_i,\nu).$ These considerations show
that 
\ba
\dist(x_1,x_2)=\inf \left(
\dist_{\p M}(x_1, z_1)+
\big(\sum_{i=1}^{p} 
t_i+\dist_{\p M}(\tilde y_i,z_{i+1})\big)\right)
\ea
where the infimum is taken over $t_i>0$,
$z_i,\tilde y_i\in \p M$, and directions  
$\xi_i,\zeta_i$ such that 
$z_{p+1}=x_2$ and the relations $(z_i,\xi_i) R_{t_i} (\tilde y_i,\zeta_i)$
are valid. }
\proofbox

\begin{theorem}
\label{boundary dist}
The boundary, $\p M$, and the broken scattering relation, $R$, determine the set ${\mathcal R}(M) \subset C(\p M)$.
\end{theorem}

\noindent
{\bf Proof.} Let $ \omega_{\p M}$ be the boundary cut locus on $M$.
{As $M \setminus \omega_{\p M}$ is dense in $M$, it is sufficient to
find ${\mathcal R}(M \setminus \omega_{\p M})$.}
Recall that, for  $x_0 \in M \setminus \omega_{\p M}$, we have $x_0=\gamma_{z_0,\nu_0}(t_0)$,
where
$t_0= \dist(x_0, \p M)<\tau_b(z_0)$ and $z_0$ is the unique boundary point closest to $x_0$.
Using the broken scattering relation $R$, we intend to determine, for any 
$w_0 \in \p M$, 
$D(z_0,t_0,w_0):=\dist(x_0,w_0)$.

Let $x(s)$ be a 
shortest path from $x_0$ to $w_0$ parametrized by the arclength. 
Denote by $w=x(s_0)$ 
the first point where $x(s)$ is in $\p M$.
Clearly, 
\beq
\label{28.7.3}
\dist(x_0,w_0) = s_0+\dist(w,w_0),\quad s_0\geq t_0.
\eeq
By \cite{Al}, the path $x([0,s_0])$ is a geodesic in $M$. We denote
$\eta= -x'(s_0)$ so that
$x_0= \gamma_{w,\eta}(s_0)$. 
As $t_0\leq \tau_b(z_0)<\tau_M(z_0)$, there is 
a family of focusing directions 
${\mathcal F}(z_0,t_0)=\{U,\xi(\,\cdotp),t(\,\cdotp)\}$
such that for $s_1=s_0$, $w_1=w$, and $\eta_1=\eta$ we have
\beq
\label{28.7.10}
(w_1,\eta_1) \, R_{s_1+t(z)}\, (z,\xi(z)),\quad z\in U.
\eeq

After these preparations we will show that 
\beq\label{eq: for D}
D(z_0,t_0,w_0)=\inf (\dist(w_0,w_1)+s_1)
\eeq
where infimum is taken over $w_1\in \p M$,
$\eta_1\in  S_{w_1}M$, and $s_1\geq t_0$
such that there is a focusing sequence
${\mathcal F}(z_0,t_0)=\{U,\xi(\,\cdotp),t(\,\cdotp)\}$
satisfying (\ref{28.7.10}). 

Formula (\ref{28.7.3}) shows that the infimum on the 
right side of (\ref{eq: for D}) is less or equal 
to $D(z_0,t_0,w_0)$. Thus to prove (\ref{eq: for D}), it is
enough to show that if $w_1,\eta_1$, and $s_1$ satisfy
(\ref{28.7.10}) then $\rho=\dist(w_0,w_1)+s_1\geq \dist(x_0,w_0)$.

{
%Let $\rho=\dist(w_0,w_1)+s_1$.  
Assume now that (\ref{28.7.10}) is valid. 
%iHOX{It seems that $r \leq t_0$ is not enough to apply L. \ref{non-inters, hausdorff}.
%Thus i slightly extended the proof. Please have a look.}
Then, for some $r(z), \tau(z), \, r(z)+\tau(z)=s_1+t(z)$, we have that
$\gamma_{z,\xi(z)}(r(z))= \gamma(\tau(z)).$

Keeping aside the
trivial case when 
the geodesics $\gamma_{z_0,\nu_0}$ and $\gamma_{w_1, \eta_1}$  
coincide, consider first the case when $\lim \sup r(z) =r >t_0$.
Denoting $\gamma_{z_0,\nu_0}(r)=x_1$, we then have
\ba
\dist(w_1,x_0)&\leq& \dist(w_1,x_1)+\dist(x_1,x_0)\\
&\leq& 
(s_1+t_0-r)+(r-t_0)\leq s_1, 
\ea
yielding 
 $\rho\geq \dist(w_0,w_1)+\dist(w_1,x_0)\geq \dist(w_0,x_0).$
 If, however,  $\lim \sup_{z\to z_0} r(z) =r\leq t_0$, we are in the situation of 
 Lemma \ref{non-inters, hausdorff}, which shows that 
 \ba
 \gamma_{z_0,\nu_0}(t_0)= \gamma(s_1),
 \ea
 yielding again that $\rho\geq  \dist(w_0,x_0).$
 %Then
%(\ref{28.7.10}) with $z=z_0$ yields that
%there is an intersection point $r\geq 0$ with
%\beq\label{eq: r is intersection}
%\gamma_{w_1,\eta_1}(s_1+t_0-r)=\gamma_{z_0,\nu_0}(r)=:x_1.
%\eeq
%we see that if $r\leq t_0$ then by Lemma \ref{non-inters, hausdorff}
%also  $r=t_0$ satisfies (\ref{eq: r is intersection}). Thus
%(\ref{eq: r is intersection}) is satisfied with some $r\geq t_0$.
%Thus (\ref{eq: for D}) is proven.}
\proofbox

As the set ${\mathcal R}(M)$ can be naturally endowed with a differential structure
and  a  Riemannian metric so that is
becomes isometric to $(M,g)$, 
see e.g.\ \cite{KaKuLa,Ku}, we have finished the proof of
Theorem \ref{main}. \proofbox

\section{Proofs for the radiative transfer equation.}

\subsection{Notations}

Let $X$ be a manifold with dimension $n$
and $\Lambda_1\subset T^*X\setminus 0$ be
a Lagrangian submanifold.
Let  $(x_1,\dots,x_n)=(x',x'',x''')$ 
of be local coordinates $X$ with $x'=(x_1,\dots,x_{d_1})$,
  $x''=(x_{d_1+1},\dots,x_{d_1+d_2})$,
 $x'''=(x_{d_1+d_2+1},\dots,x_{n})$,  and
$\phi(x,\theta)$, $\theta\in \R^N$ be a non-degenerate
phase function that parametrizes $\Lambda_1$.
We say that distribution $u\in {\cal D}'(X)$ 
is a Lagrangian distribution associated with $\Lambda_1$ and denote
$u\in I^m(X;\Lambda_1)$, if it can can locally be represented
as
\ba
u(x)=\int_{\R^N}e^{i\phi(x,\theta)}a(x,\theta)\,d\theta,
\ea
where $a(x,\theta)\in S^{m+n/4-N/2}(X\times \R^N\setminus 0)$,
see  \cite{GU1,H4,MU1}.

Let $S_1\subset X$ be a submanifold of codimension $d_1$. 
We denote its conormal bundle by
$N^*S=\{(x,\xi)\in T^*X\setminus 0:\ x\in S,\ \xi\perp T_xS\}$.
If $S_1=\{x'=0\}$ in local coordinates, $\Lambda_1=N^*S_1$
and $u\in I^m(X;\Lambda_1)$,
then locally
\ba
u(x)=\int_{\R^{d_1}}e^{i x'\cdotp \theta}a(x,\theta')\,d\theta',
\quad
a(x,\theta')\in S^{\mu}(X\times \R^{d_1}\setminus 0)
\ea
where $\mu=m-d_1/2+n/4$. We denote
 $I^m(X;S_1)=I^\mu(X;N^*S_1)$ and say that
$I^\mu(X;S_1)$ are the conormal distributions in space $X$
associated with submanifold $S_1$.

Also, we denote by $I^{p,l}(X;\Lambda_1,\Lambda_2)$
the distributions $u$ in ${\cal D}'(X)$ associated to two
cleanly intersecting Lagrangian manifolds $\Lambda_1,\Lambda_2\subset T^*X$, see \cite {GU1,MU1}.
Let $S_1$ and $S_2$ be submanifolds of $M$ of codimensions
$d_1$ and $d_1+d_2$,
$S_2\subset S_1$. If in local coordinates
$S_1=\{x'=0\}$, $S_2=\{x'=x''=0\}$,
and $\Lambda_1=N^*S_1,$  $\Lambda_2=N^*S_2$, then
the distribution $u\in I^{p,l}(X;\Lambda_1,\Lambda_2)$
can be locally represented as
\ba
u(x)=\int_{\R^{d_1+d_2}}e^{i (x'\cdotp \theta'+
x''\cdotp \theta'')}a(x,\theta',\theta'')\,d\theta'd\theta'',
\ea
where $
a(x,\theta',\theta'')$ belongs to a product type symbol class
$S^{\mu',\mu''}(X\times (\R^{d_1}\setminus 0)\times \R^{d_2})$
containing symbols $a\in C^{\infty}$ that satisfy
\ba
|\p_x^\gamma\p_{\theta'}^\alpha \p_{\theta''}^\beta a(x,\theta',\theta'')|
\leq C_{\alpha\beta\gamma K}(1+|\theta'|+|\theta''|)^{\mu-|\alpha}|
(1+|\theta''|)^{\mu'-|\beta|}
\ea
for all $x\in K$, multi-indexes  $\alpha,\beta,\gamma$, and
compact sets $K \subset X$. Above,
$\mu=p+l-d_1/2+n/4$ and
$\mu'=-l-d_2/2$.

By \cite{GU1,MU1}, microlocally away from $\Lambda_1\cap \Lambda_2$,
\ba
I^{p,l}(\Lambda_0,\Lambda_1)\subset I^{p+l}(\Lambda_0\setminus \Lambda_1)\quad
\hbox{and}\quad I^{p,l}(\Lambda_0,\Lambda_1)\subset I^{p}(\Lambda_1\setminus \Lambda_0).
\ea
Thus the principal symbol of $u\in I^{p,l}(\Lambda_0,\Lambda_1)$
is well defined on $\Lambda_0\setminus \Lambda_1$ and  $\Lambda_1\setminus \Lambda_0$.

\subsection{Born series}

In the sequel, we denote the distance on $(N,g)$ by $d(x,y)=\dist(x,y)$.
Let $\gamma_{x,\xi}(t)$ be the geodesic on $(N,g)$
with initial point $x$ and initial direction $\xi\in S_{x_0}N$.
Denote 
\ba% \label{notat a}
& &\gamma_{x,\xi}=\{\gamma_{x,\xi}(t)\in N:\
t\in \R\},\\ \nonumber
& &\eta_{x,\xi}=\{(\gamma_{x,\xi}(t),\dot \gamma_{x,\xi}(t))\in SN:\
t\in \R\},\\
& & \nonumber
\eta^+_{x,\xi}=\{(\gamma_{x,\xi}(t),\dot \gamma_{x,\xi}(t))\in SN:\
t\in \R_+\}.
\ea
The measurement operator $A$ can be extended to distributions $w$
supported in $SU$. 
In the following we consider   $u$ corresponding
to $w_0(x,\xi)=\delta_{(x_0,\xi_0)}(x,\xi)$, $x_0\in U$.
We assume that $\gamma_{x_0,\xi_0}(\R_+)$ intersect the strictly convex manifold
$M\subset N$.
To analyze the corresponding solution, let us denote
the specific geodesic on which the leading
order singularities propagate by
$\gamma_0=\gamma_{x_0,\xi_0}$. 
Also, we denote the corresponding spray in $SN$ by
$\eta_0=\eta_{x_0,\xi_0}$.

Let  $u_0(x,\xi,t)$ be the solution of 
the equation (\ref{a1}) with $S$ being zero, that is, $Hu_0+\sigma u_0=0$, 
$u_0|_{t=0}=w_0$. Then $u_0(t)=c_0(x)\delta_{\eta_0(t)}(x,\xi)$,
where $c_0(x)$ is a non-vanishing smooth function. 
To simplify notations,
we consider the equation for all $t\in \R$, obtaining
\ba
u_0(t,x,\xi)=c_0(x)\delta_{\eta_0(t)}(x,\xi),\quad (t,x,\xi)\in \R\times SN.
\ea

In the following we analyze the higher order
terms in the Born series, that is,
\ba
& &u_j=QSu_{j-1},\quad j\geq 1,
\ea
where $Q$ is defined by $v=QF$ where 
\ba
& &Hv+\sigma v=F\quad\hbox{in }\R_+\times SN,\quad\quad v|_{t=0}=0.
\ea

We note that there are $C_1,C_2>0$ so 
the solutions $u^w$ of equation (\ref{a1})
satisfy  
\beq\label{exponetial growth}
|u^w(x,\xi,t)|\leq C_1e^{C_2t}\|w\|_{L^\infty(SN)}.
\eeq

To analyze the singularities of $u$, let
us take the Laplace transform $\L$
in time $t$ and 
consider $\hat u(k,x,\xi)=(\L u(\cdotp,x,\xi))(k)$.
By (\ref{exponetial growth}) the Laplace
transform is well defined for $k\in \C$, Re\,$k>C_2$. 
In the following, we consider $k$ first as as a parameter,
and denote $\hat u(x,\xi)=\hat u(x,\xi,k)$. Then
\ba
& &(k+\hat H)\hat u+\sigma \hat u-S\hat u=w_0\quad\hbox{in }(x,\xi)\in SN,
\ea
where $w_0(x,\xi)=\delta_{(x_0,\xi_0)}(x,\xi)$ and
\ba
\hat H u(x,\xi)=-\xi^j\frac {\p u}{\p x^j}-\xi^l\xi^j 
\Gamma_{lj}^m(x)\frac {\p u}{\p \xi^m}.
\ea
The operator $\hat H+k:C^\infty(SN)\to C^\infty(SN)$
has  $\hat Q_k$ a parametrix, see \cite{H4,MU1},
that  satisfies $\L(QF)(k)=\hat Q_k (\L F(k))$. 
Also, we denote
$\hat u(k)=\hat u_0(k)+\hat u_{sc}(k)$, where
$\hat u_{sc}(k)=\hat u_1(k)+\hat u_2(k)+\dots$.

Consider now a Born iteration
 starting at
a general $w_0(k)$.
Since the coefficients
of $\hat H$ are smooth functions and the kernel of $S$ is a smooth
compactly supported function, we that 
for any $s\geq 0$ there there is $C_3=C_3(s)>0$ such that
for  Re\,$k>C_3$  
the Born series
\beq
\hat w(k)=\sum_{j=0}^\infty (\hat Q_k S)^{j-1}w_0(k) 
\label{Born series II}
\eeq 
 converges in Sobolev space $H^s_{loc}(SN)$  
when $w_0(k)\in H^s_{loc}(SN)$.

\subsection{Properties of the compositions of the operators $S$ and $\hat Q_k$}

\begin{lemma}
We can write $S=S_1S_2$,
\ba
S_jf(x,\xi)=\int_{S^{n-1}}K_j(x,\xi,\xi')f(x,\xi')\,dS(\xi'),\quad
j=1,2
\ea
where $K_j(x,\xi,\xi')\in C^\infty_0(SN\times SN)$. 
\end{lemma}

\noindent
{\bf Proof.}
Interpreting $x$ as a parameter, we define 
$K_x:L^2(S^{n-1})\to L^2(S^{n-1})$ by 
\ba
K_xf(\xi)=\int_{S^{n-1}}K(x,\xi,\xi')f(\xi')\,dS(\xi').
\ea

As the kernel $K(x,\xi,\xi')$ is smooth,
we see that for all $\alpha\in \N^n$ and $l,m\in \N$
there is a constant $c_{\alpha l m}$ such that
\beq\label{suppi}
\sup_{x\in M}\|\p_x^\alpha(1-\Delta_{\xi})^mK(x,\xi,\xi')\|_{C^l(
S^{n-1}\times S^{n-1})}
<c_{\alpha l m},
\eeq
where $\Delta_\xi$ is the Laplace-Beltrami operator
of the $(n-1)$-sphere $S^{n-1}$. Let $a_m>0$ be numbers such that
$0<a_m<e^{-m}\min(1,c_{\alpha l m}^{-1})$ for all
$\max(|\alpha|,l)\leq m$.
 Then the
operator
\ba
B=\sum_{m=0}^\infty a_m(1-\Delta_{\xi})^m
\ea
defines an unbounded non-negative selfadjoint 
operator $B:L^2(S^{n-1})\to L^2(S^{n-1})$ having an
inverse $J=B^{-1}$ that can be extended to a smoothing
operator ${\cal D}'(S^{n-1})\to C^\infty(S^{n-1})$.
Moreover, by (\ref{suppi}) we see that for any $x$
the operator $L_x=BK_x$ defines 
a smoothing operator  ${\cal D}'(S^{n-1})\to C^\infty(S^{n-1})$
and its Schwartz kernel $L_x(\xi,\xi')$ is a $C^\infty$-smooth
in all variables $(x,\xi,\xi')$. Thus we prove the assertion
by defining $K_2(x,\xi,\xi')=L_x(\xi,\xi')$
and $K_1(x,\xi,\xi')=J(\xi,\xi')$, where
$J(\xi,\xi')$ is the Schwartz kernel of $J$.
\proofbox

The Born series iteration can be written  as
\ba
\hat u_j(k)=\hat Q_kS_1A^{j-1}S_2\hat u_0(k)
\ea
where $ A=S_2\hat Q_kS_1$. To analyze the operator $A$
we consider first the case where
$K(x,\xi,\xi')$ would be the constant $1$.
Denote by $S^c$ the  operator corresponding to
a constant scattering kernel $K(x,\xi,\xi')=1$. 
For this purpose, we introduce
operators $T=\pi_*:L^2(SN)\to L^2(N)$ and  $T^*=\pi^*:L^2(N)\to L^2(SN)$,
that is,
\ba
Tu(x)=c_n^{-1}\int_{S_xN}u(x,\xi)dV_g(\xi),\quad\quad T^*v(x,\xi)=v(x),
\ea
where $ c_n=\hbox{vol}(S^{n-1})$ and $V_g$ is the volume on $S_xN$.

\begin{lemma}\label{lem: A compositions}
Let  $Z=SN\times SN$,
$L_0=\{(x,\xi,y,\eta)\in Z:\ x=y\}$, and $\Sigma_0=N^*L_0$. 
The Schwartz kernels of $A^c$ and $A$ satisfy
\beq\label{A^c-condition}
& &A^c(x,\xi,y,\eta)\in I^{-1}(Z;L_0)=I^r(Z;\Sigma_0),\\
\label{A-condition}
& &A(x,\xi,y,\eta)\in I^{\rho}(Z;\Sigma_0)
\eeq
where $r=-(n+1)/2$,
%$r=-1+\frac n2-\frac {2(2n-1)}4=-(n+1)/2$,
$\rho=r+\e$,
and $\e>0$.
\end{lemma}

\noindent
{\bf Proof.}
Clearly, $TT^*=I$ and $S^c=T^*T$.
Thus we have $S^c=S^c_1S^c_2$ where
$S^c_1=S^c_2=S$.
In the local 
coordinates $S^c$ has the Schwartz kernel
\ba
S^c(x,\xi,x',\xi')=\delta(x-x')\in I^0(Z;L_0)=I^{m_1}(Z;\Sigma_0),
\ea
where $m_1=(1-n)/2$.
To analyze $A=S_2\hat Q_kS_1$, we first consider the operator
\ba
A^c=S^c_2\hat Q_kS^c_1=T^*T \hat Q_kT^*T.
\ea
Denote  $\tilde Q_k=T \hat Q_kT^*:L^2(N)\to L^2(N)$ 
and let $v\in C^\infty_0(N)$. Then
\ba
(\hat Q_k T^*v)(x,\xi)=\int_{-\infty}^0  
h(x,\xi,s,k)
v(\gamma_{x,\xi}(s))\,ds
\ea
where $h(s,x,\xi,k)$ is the solution of the 
differential equation
\beq\label{ODY}
& &\p_s h(s,x,\xi,k)+(k+\sigma(\gamma_{x,\xi}(s))) h(s,x,\xi,k)=0,\\
& &h(s,x,\xi,k)|_{s=0}=1.\nonumber
\eeq
Note that
\beq\label{k dependency}
h(x,\xi,s,k)=e^{-ks}h(x,\xi,s,0).
\eeq
Thus, using the assumption that the manifold $N$ is simple, we have
\beq\label{tilde Q is PsiDO}\quad\quad
(T\hat Q_kT^* v)(x)&=&
\int_{S^{n-1}}
\int_{-\infty}^0  
h(s,x,\xi,k)
v(\gamma_{x,\xi}(s))\,dsdV_g(\xi)\\
\nonumber
&=&\int_{N}
[h(s(x,y),x,\xi(x,y),k)j(x,y)] v(y)\, dV_g(y),
\eeq
where $s(x,y)\in (-\infty,0]$ and $\xi(x,y)\in S_xN$ 
are defined by $\exp_{x}^{-1}(y)=s(x,y)\xi(x,y)$,
and $j(x,y)=\det(d\exp_x|_y)^{-1}$ is the Jacobian determinant
where $d\exp_x|_y$ is the differential of the map
$\exp_x$ evaluated at $y$. Since $(N,g)$ is simple,
the kernel $b(x,y):=h(s(x,y),x,\xi(x,y),k)j(x,y)$ is smooth
outside the diagonal and 
behaves near the diagonal as 
\ba
b(x,y)\sim e^{-kd(x,y)}d(x,y)^{1-n}.
\ea 
Using (\ref{tilde Q is PsiDO}) we see
that $\tilde Q_k$ is a pseudodifferential operator
of order $(-1)$ (for a similar argument see \cite{SU1}).

The Schwartz kernel  $\tilde Q_k(x,x')\in I^{-1}(N\times N;
\diag(N\times N))$ of $\tilde Q$ can be written
as 
\ba
\tilde Q_k(x,x')=\int_{\R^n}e^{i(x-x')\cdotp \theta}a(x,x',\theta)d\theta,
\quad a\in S^{-1} (N\times N\times \R^n\setminus 0).
\ea
The same expression defines a function $\tilde Q_k(x,\xi,x',\xi')
:=\tilde Q_k(x,x')\in I^{-1}(SN\times SN;
L_0)$. This function is the Schwartz kernel of $A^c=T^*\tilde Q_kT$ and thus 
we see that the first part of the assertion, the formula
(\ref{A^c-condition}) is  satisfied.

Next we consider the Schwartz kernel of $A$, that is,
$A(x,\xi,y,\eta)$. It can be written as a product
\ba
A(x,\xi,y,\eta)=A^c(x,\xi,y,\eta)J(x,\xi,y,\eta)
\ea
where (using the Riemannian normal coordinates at $x$)
\ba
J(x,\xi,y,\eta)=K_2(x,\xi,\frac {y-x}{|y-x|})K_1(y,\frac {x-y}{|x-y|},\eta).
\ea

Now $J_1(x,y,z):=K_1(x,\xi,z/|z|)$ and $J_2(x,y,z):=K_2(x,z/|z|,\xi)$ 
are homogeneous functions if degree zero in $z$, and we 
see that \cite[formula (1.2)]{GU1} 
\ba
K_2(x,\xi,\frac {y-x}{|y-x|}), K_1(y,\frac {x-y}{|x-y|},\eta)
\in I^{-n}(Z;L_0).
\ea
Now we can write $A$ as the product of $K_1$, $K_2$, and $A^c$.
To analyze this product, we need the following lemma
extending results of \cite{GU1} for less regular conormal distributions.

\begin{lemma}\label{lem: products}
Let $Z$ be a manifold of dimension $d$ and $L_0$ be a 
submanifold with codimension $n$.  Assume that
$A\in I^{-d}(Z;L_0)$ and $B\in I^{\mu}(Z;L_0)$, $\mu<0$.
Then the pointwise product $A B \in I^{\mu+\e}(Z;L_0)$ for any $\e>0$.
\end{lemma}

\noindent
{\bf Proof.} Let $(z',z'')$ be local coordinates of $X$ such that
$L_0=\{z'=0\}$. Then 
\ba
A(z)=\int_{\R^d}e^{iz'\cdotp \theta}a(z,\theta)\,d\theta,\quad
B(z)=\int_{\R^d}e^{iz'\cdotp \theta}b(z,\theta)\,d\theta,\quad
\ea
where $a(z,\theta)\in S ^{-d}(X\times \R^d\setminus 0)$ and  
$b(z,\theta)\in S ^{\mu}(X\times \R^d\setminus 0)$. The 
symbol $c(z,\theta)$  of the product $A(z)B(z)$ is given by the convolution
\ba
c(z,\theta)=\int_{\R^d}a(z,\theta-\tilde \theta)\,
b(z,\tilde \theta)\,d\tilde \theta,
\ea
and a simple computations shows that
\ba
|c(z,\theta)|\leq C
\int_{\R^d} (1+|\theta-\tilde \theta|)^{\mu}
(1+|\tilde \theta|)^{-d}\,d\theta\leq C'(1+|\theta|)^{\mu+\e},
\ea
with $\e>0$.
Indeed, decomposing the domain of integration as
$\R^d=B(0,\frac 12|\theta|)\cup B(\theta,\frac 12|\theta|)\cup
\big(\R^d\setminus (B(0,\frac 12|\theta|)
\cup B(\theta,\frac 12|\theta|))\big)$,
we see that 
\ba
|c(z,\theta)|&\leq& 
C_1 |\theta|^{\mu}\log |\theta|+
C_2  | \theta|^{-d} |\theta|^{d+\mu}(1+\delta_{\mu,-d}\log|\theta|)+
C_3 |\theta|^{\mu}\\
& \leq& C'(1+|\theta|)^{\mu+\e},
\ea
where $|\theta|>1$ and $\delta_{\mu,-d}$ is one if $\mu=-d$ and
zero otherwise.
 The derivatives of $c(z,\theta)$ can be estimated in similar
way, and we obtain that  $c(z,\theta)\in S ^{\mu+\e}(X\times \R^d\setminus 0)$.
\proofbox

Lemma \ref{lem: products} for the product of 
 $K_1$, $K_2$, and $A^c$ implies (\ref{A-condition}).
This proves Lemma \ref{lem: A compositions}.\proofbox

The previous result says, roughly speaking,  that
 $A$ is like a $\Psi$DO of order $(-1)$ operating
in $(x,y)$-variables when $\xi$ and $\eta$ are considered
as parameters. 

Next we consider powers of $A$. 
Next, $\Sigma_0'$ denotes the canonical relation
corresponding to the Lagrangian manifold $\Sigma_0$.
We see that
$\Sigma_0'\times \Sigma_0'$ intersects 
cleanly $T^*SN\times \diag(T^*SN\times T^*SN)\times T^*SN$
with the excess $d=(n-1)$. 
Thus using \cite[Thm VIII.5.2]{Tr2},
we see that
\ba
A^2=A\circ A\in I^{-2\rho+d/2}  (Z;\Sigma_0)=
I^{\rho_2}  (Z;\Sigma_0),
\ea
where $\rho_2=-(n+3)/2+2\e$ with any $\e>0$.
Iterating operator $A$,
we see that
\ba
%\label{eq: b-kaava} \quad\quad\quad
A^j\in I^{\rho_j}  (Z;\Sigma_0)=
I^{-1-j+\e}  (Z;L_0),\ \ \rho_j=-\frac {n+1}2-j+\e,\ \e>0.\hspace{-1cm}
\ea

\subsection{Singularities of the terms in the Born series}

In the following, 
let
$\Lambda_0=N^*Y_0$ and
$\Lambda_1=N^*(Y_1)$, 
where
\ba
Y_0&=&\{(\gamma_0(t),\dot \gamma_0(t))\in SN:\ t\in \R\},\\
Y_1&=&S(\gamma_0)
=\{(x,\xi)\in SN:\ x\in \gamma_0(\R)\}.
\ea
Moreover, let $P=P(x,\xi,D_x,D_{\xi})=\hat H+k$,
\ba
\char(P)=\{(x,\xi,\tilde x,\tilde \xi)\in T^*(SN):\
\xi^i\tilde x_i+\xi^i\xi^j 
\Gamma_{ij}^k(x)\tilde \xi_k=0\},
\ea
and let $\Xi(x,\xi,\tilde x,\tilde \xi)$ be the bicharacteristic of 
$P(x,\xi,D_x,D_{\xi})$ (i.e.\ integral curve of the Hamilton vector 
field in $T^*(SN)\setminus 0$) starting from $(x,\xi,\tilde x,\tilde \xi)\in T^*(SN)$.
Then the  flow-out canonical relation generated by 
 $\char(P)$ is
\ba
& &\Lambda_P'=\{(x,\xi,\tilde x,\tilde \xi;y,\zeta,\tilde y,\tilde \zeta)
\in (T^*(SN)\setminus 0)\times
 (T^*(SN)\setminus 0):\\
& &\quad \quad \quad \quad \quad \quad \quad \quad (x,\xi,\tilde x,\tilde \xi)\in 
\char(P),\ \ 
(y,\zeta,\tilde y,\tilde \zeta)\in \Xi(x,\xi,\tilde x,\tilde \xi)\}.
\ea
The flow-out of $\Lambda_1$ in $\char(P)$ is the
Lagrangian manifold $\Lambda_2\subset T^SN\setminus 0$
 satisfying $\Lambda_2'=\Lambda_P'\circ \Lambda_0'$.

\begin{lemma} We have 
\ba
%\label{0-term}
\hat u_0(x,\xi,k)=c_0(x,k)\delta_{\eta_0}(x,\xi)
\in I^{r_0}(SN;\Lambda_0),
\ea
where $c_0(x,k)$ is a smooth non-vanishing function and
%$r_0=0+\frac{2n-2}2-\frac{2n-1}4=
%\frac{2n-3}4. $
$r_0=(2n-3)/4. $
For $j\geq 1$, 
\beq\label{Born k-term}\quad\quad
\hat u_j(k)\in I^{r_j,-\frac 12}(SN;\Lambda_1,\Lambda_2),\ \ 
r_j=-\frac n2-j+\e\delta_{j\geq 2},\ \e>0,\hspace{-1cm}
\eeq
where $\delta_{j\geq 2}$ is one if $j\geq 2$ and zero otherwise.
\end{lemma}

\noindent
{\bf Proof.}
For the  zeroth term in the Born series 
the claim is true by definition.
Next we analyze the higher order terms.
Clearly, 
\ba
S_2\hat u_0(x,\xi,k)=
K_2(x,\xi,\eta(x))(S^c\hat u_0)(x,\xi,k),
\ea
 where $\eta(x)\in S_xN$ defines a smooth vector
field such that if $x=\gamma_0(s)$ then $\eta(x)=\dot \gamma_0(s)$.
A simple computation shows that
$\Lambda_0'\times \Sigma_0'$ intersects
$\hbox{diag}(T^*SN\times T^*SN)\times (T^*SN)$ transversally.
Now
$S_2\in I^{0}(SN\times SN;L_0)=I^{m_1}(SN\times SN;\Sigma_0)$, where
%$m_1=0+\frac n2-\frac{2(2n-1)}4=\frac {1-n}2$ 
$m_1=(1-n)/2$ 
and
 by \cite[Thm 25.2.3]{H4} that 
 $S_2$ can be considered as a continuous operator
\ba%\label{first S operation 3}
S_2: I^{r_0}(SN;\Lambda_0)\to I^{s}(SN;\Lambda_1),
\ea
where $s=r_0+m_1$
and $\Lambda_1'=\Lambda_0'\circ \Sigma_0'$.
A simple computation shows that
$\Lambda_1'\circ\Sigma_0'=\Lambda_1'$, and that
$\Lambda_1'\times\Sigma_0'$ intersects
$\hbox{diag}(T^*SN\times T^*SN)\times (T^*SN)$ cleanly
with excess $e=(n-1)$.
Thus we have by \cite[Thm 25.2.3]{H4} that 
\ba%\label{eq: A terms}
A^jS_2\hat u_0(k)\in  I^{\rho_j+m_1+e/2}(SN;\Lambda_1).
\ea
Again, as
$\Lambda_1'\circ\Sigma_0'=\Lambda_1'$, and
$\Lambda_1'\times\Sigma_0'$ intersects
$\hbox{diag}(T^*SN\times T^*SN)\times (T^*SN)$ cleanly
with excess $e$, we see that
since $S_1\in I^{m_1}(Z;\Sigma_0)$, 
\beq\label{eq: A terms2}\quad
S_1A^jS_2\hat u_0(k)\in  I^{\rho_j+2(m_1+e/2)}(SN;\Lambda_1)
= I^{\rho_j}(SN;\Lambda_1).
\eeq
To analyze $\hat u_j(k)=\hat Q_kS_1A^{j-1}S_2\hat u_0(k)$,
we observe that the operator $P=\hat H+ik$ is a first order operator
of real principal type. As $\hat Q_k$
is its parametrix, it follows from 
\cite{MU1} that the Schwartz kernel
\beq\label{I class of hat Q_k}
\hat Q_k\in I^{\frac 12-1,-\frac 12}(Z;\Delta_{T^*Z},\Lambda_{P}),
\eeq
where $\Delta_{T^*Z}'$ is the diagonal of $T^*Z\times T^*Z$ and 
 $\Lambda_{P}'\subset T^*(Z)$ 
is the flow-out canonical relation generated by $\char(P)$.
Now $N^*Y_1$ intersects $\char(P)$ transversally.
Hence we obtain (\ref{Born k-term}) by  
\cite[Prop. 2.1]{GU1}.
\hfill$\Box$

\subsection{Principal symbol of the singularity}

For any $s>0$ there is $j_0$ such that
$\hat u_{j_0}(k)\in H^s_{loc}(SN)$. Using the convergence
of the Born series (\ref{Born series II}), we 
see that the series 
$\hat u_{j_0}(k)+\hat u_{j_0+1}(k)+\hat u_{j_0+2}(k)+\dots$
converges in $H^s_{loc}(SN)$.

Next we consider how to find the geodesic $\gamma_0$ in $U$. 
To this end we observe using (\ref{eq: A terms2})
that  $T\hat u(k)=
T\hat u_0(k)+T\hat u_{sc}(k)\in I^0(N;\gamma_0)$ 
and $T\hat u_0(k)\in I^0(N;\gamma_0)$ have the same non-vanishing 
principal symbol. Thus $T\hat u(k)$ in $U$ determines  $U\cap \gamma_0$.

Moreover, the above convergence of the Born series in Sobolev spaces
and  (\ref{Born k-term}) yield that
$\hat u_1(k)$ and $\hat u_{sc}(k)=\hat u_1(k)+\hat u_2(k)+\dots$ are both
elements in $I^{r_1,-\frac 12}(SN;\Lambda_1,\Lambda_2)$
and they have the same principal symbol on $\Lambda_2\setminus\Lambda_1$.
Motivated by this,  we consider next $\hat u_1(k)$.  

Using the above notations, we see that
\ba
S\hat u_0(x,\xi,k)=
S(x,\xi,\eta(x))h(d(x,x_0),x_0,\xi_0,k)c_1(x)\delta_{\gamma_0}(x)
\in I^0(SN;Y_1),
\ea
where $c_1(x)$ is a smooth non-vanishing function.
Moreover, the operator $\hat Q_k$ has the Schwartz kernel
(\ref{I class of hat Q_k}) that 
away from the diagonal has the form
\ba
\hat Q_k(x,\xi,x',\xi')= 
h(d(x,x'),x',\xi',k) \delta_{\eta^+_{x',\xi'}}(x,\xi),
\ea
where $h$ is defined in (\ref{ODY}).
Thus, in $(x,\xi,x',\xi')\in Z\setminus L_0$, the kernel of $\hat Q_k$  has the form
\ba
& &\hat Q_k(x,\xi,x',\xi')=\\
& &
\int_{\R^N} e^{i\psi(x,\xi,x',\xi',\theta)}
[h(d(x,x'),x',\xi',k)
q(x,\xi,\theta)]\,d\theta \quad \hbox{mod }C^\infty(Z)
\ea
where $\psi(x,\xi,x',\xi',\theta)$ is a non-degenerate phase function
parameterizing the Lagrangian $\Lambda_P$
and  $q(x,\xi,\theta)\in S^{r_1-1/2+(4n-2)/4-N/2} (\R^n\times \R^{n-1}\times
\R^N\setminus 0)$ has a non-vanishing principal symbol.

Let us use in $SN\setminus \eta_0$ local coordinates 
${\mathcal S}:(x,\xi)\mapsto (s_j(x,\xi))_{j=1}^{2n-1}$ having the property
that if $\gamma_{x,\xi}(\R_-)$ intersects the geodesic $\gamma_0(\R_+)$
then $s_1=s_1(x,\xi)$ is the unique value such
that 
\ba
\gamma_{x,\xi}(\R_-)\cap \gamma_0(\R_+)=\gamma_0(s_1),
\ea
and $s_2(x,\xi)=d(\gamma_0(s_1(x,\xi)),x)$.
By \cite[Prop. 2.1]{GU1}, \ba
\hat u_1(k)=\hat Q_kS\hat u_0(k)
\in I^{r_1,-\frac 12}(SN;\Lambda_1,\Lambda_2)
\ea 
and $\hat u_1(x,\xi,k)$ in $(x,\xi)\in SN\setminus \eta_0$ has in the  
above local coordinates the form
\ba
& &\hat u_1(x,\xi,k)=
\int_{\R^N}e^{i\phi(x,\xi,\theta)}
[a(x,\xi,k)p(x,\xi,\theta)]\,d\theta\quad \hbox{mod }C^\infty(SN),\\
& &a(x,\xi,k)= h\big(s_1(x,\xi),x_0,\xi_0,k\big)\,h\big(s_2(x,\xi),
\gamma_0(s_1(x,\xi)),\zeta(x,\xi),k\big)
\ea
where $\phi(x,\xi,\theta)$ is a non-generate phase function parametrizing
the Lagrangian manifold
$\Lambda_2$, 
 $\zeta(x,\xi)=-\dot\gamma_{x,\xi}(-s_2(x,\xi))$
is the direction of $x$ from $\gamma_0(s_1)$              
and $p(x,\xi,\theta)$ is a symbol with a non-vanishing principal symbol.
Note that on $\Lambda_2\setminus \Lambda_1$ the
principal symbol of $a(x,\xi,k)p(x,\xi,\theta)$ is non-vanishing
on the conormal bundle of the submanifold
\ba
K=\{(x,\xi)\in SN:\ \gamma_{x,\xi}(\R_-)\cap \gamma_0(\R_+)\cap M^{\rm int}\not=\emptyset\}.
\ea
By (\ref{k dependency}),
\beq\label{eq: asymptotics}
a(x,\xi,k)
=e^{-k(s_1+s_2)}\,
S(\gamma_0(s_1),\zeta,\dot \gamma_0(s_1))\,b_0(x,\xi),
\eeq
where $s_1=s_1(x,\xi)$, $s_2=s_2(x,\xi)$,
$\zeta=\zeta(x,\xi)$,
 and $b_0(x,\xi)$ is non-vanishing and independent of $k$. 

Now we are ready prove unique solvability of the inverse problem.
\smallskip

\noindent
{\bf Proof of Theorem \ref{transport-main}}.
First we note that have found already the 
set $\gamma_0\cap U$. 
Thus we know
the set $W:=SN\setminus  (SM\cup \eta_0)$.
By observing the singularities of $\hat u(k)$ at $W$, 
we can find the conormal bundle of the manifold $K\cap U$.
 Thus
by observing $\hat u(k)$ at $W$
 we can find all  points $(x,\xi)\in W$ such that
there is a broken geodesic from $(x_0,\xi_0)$
to $(x,\xi)$ with a breaking point in $M^{\rm int}$. Moreover, 
we can find
the principal symbol 
of $\hat u(k)$ on $N^*K \cap W$ 
in some local coordinates. 
 By (\ref{eq: asymptotics}), 
observing the asymptotics of the principal symbol
on  $N^*K \cap W$
when $k\to \infty$,
we can find the function $d(x_0,\gamma_0(s_1))+d(\gamma_0(s_1),x)$,
$s_1=s_1(x,\xi)$
on $(x,\xi)\in W$. 
Here $\gamma_0(s_1)\in M^{int}$
is the point at which the broken geodesic
from $(x_0,\xi_0)$
to $(x,\xi)$ breaks, that is,
the broken geodesic changes its direction.

Using the continuity of the geodesic flow, we can
find all $(x,\xi)\in SN\setminus SM$ that are in
the broken scattering relation with $(x_0,\xi_0)$ and moreover,
in such case we can find 
the broken geodesic distance
$d(x_0,\gamma_0(s_1))+d(\gamma_0(s_1),x)$.
This proves the result and even more: 
The singularities of the Schwartz kernel of the operator $A$
determine the broken scattering relation $R$.
\proofbox

\noindent
{\bf Acknowledgements:} M.\ Lassas was partially supported 
by the Academy of Finland Center of
Excellence programme 213476. 
G.\ Uhlmann was partially supported by FRG grant DMS 0554575
and a Walker Family Endowed Professorship.

\end{document}